\documentclass[11pt]{amsart}

\usepackage{geometry}
\usepackage{amsmath}

\usepackage{xcolor}
\usepackage{amsmath}
\usepackage{amssymb}
\usepackage{amsfonts}
\usepackage{graphicx}
\usepackage{booktabs}
\usepackage{algorithm}
\usepackage{algorithmicx}
\usepackage{algpseudocode}
\usepackage{enumitem}
\usepackage{subcaption}
\usepackage{microtype}

\usepackage{amsopn}

\usepackage{graphicx}  
\usepackage{epstopdf}%This is to transfer .eps figure to .pdf figure; please compile your paper using PDFLeTex or PDFTeXify.
\usepackage{epstopdf}
\usepackage{verbatim}
\usepackage{mathrsfs}
\usepackage{mathtools}
\usepackage{pstricks}
\usepackage{relsize}
\usepackage{tikz}
\usetikzlibrary{matrix}
\usepackage{subcaption}
\usepackage{pgfplots}
\usepackage{fixltx2e}
\usepackage{enumitem}
\usepackage{ upgreek }
\usepackage[colorlinks=true]{hyperref}
\hypersetup{urlcolor=blue, linkcolor=blue, citecolor=red}
\usepackage{bm}
\usepackage{appendix}

\usepackage{xcolor}

\usepackage{soul}
\usepackage{color}

\usepackage[abbrev]{amsrefs}
\usepackage{amssymb}

\usepackage[all,cmtip]{xy}

\usepackage{tikz}
\usepackage{tikz-cd}

\numberwithin{equation}{section}

\theoremstyle{plain} %% This is the default, anyway
\newtheorem{theorem}{Theorem}[section]
\newtheorem{lemma}[theorem]{Lemma}

\newtheorem{proposition}[theorem]{Proposition}
\newtheorem{remark}[theorem]{Remark}
\newtheorem{definition}[theorem]{Definition}

\definecolor{ForestGreen}{RGB}{34,139,34}
\definecolor{ao(english)}{rgb}{0.0, 0.5, 0.0}

\usepackage{lineno}
\begin{document}

\title[Adaptive Approximation of High-Dimensional HJB Equations]
{Residual-Driven Basis and Data Selection for Polynomial
Approximation of High-Dimensional Hamilton--Jacobi--Bellman Equations}

\author[Y. Huang]{Yuyang Huang\textsuperscript{1}}
\author[D. Kalise]{Dante Kalise\textsuperscript{2}}
\author[Y. Zhu]{Yuhua Zhu\textsuperscript{3}}
\thanks{\noindent\textsuperscript{1}\,Mathematical Institute, University of Oxford, Oxford OX2 6GG, UK.
\texttt{yuyang.huang@maths.ox.ac.uk}.
\newline
\textsuperscript{2}Department of Mathematics,
Imperial College London, South Kensington Campus,
London SW7 2AZ, UK.
\texttt{d.kalise-balza@imperial.ac.uk}.
\newline
\textsuperscript{3}Department of Statistics and Data Science,
University of California, Los Angeles, USA.
\texttt{yuhuazhu@ucla.edu}.%
}
\date{}

% Insert your existing abstract here:
%
% \begin{abstract}
% Your existing abstract.
% \end{abstract}
%
% Insert any existing keywords and MSC classifications here:
% \keywords{...}
% \subjclass[2020]{...}
 
\begin{abstract}
We develop greedy strategies for selecting global polynomial basis functions and training data in the sparse approximation of high-dimensional partial differential equations, with emphasis on Hamilton--Jacobi--Bellman (HJB) equations arising in nonlinear optimal control. Inspired by the goal-oriented dual-weighted residual framework, we derive a closed-form, mesh-free indicator that predicts, to first order, the reduction of the PDE residual obtained by adding a single candidate basis function, without solving an adjoint problem. In the data-driven regime, we jointly select basis functions and gradient-augmented training samples by maximizing the predicted decrease of a supervised loss. Numerical experiments on benchmark PDEs and nonlinear optimal control problems in up to $35$ dimensions show accuracy comparable to non-adaptive approaches with substantially fewer basis functions and training samples.
\end{abstract}


\maketitle

\section{Introduction}
The numerical solution of partial differential equations (PDEs) posed on high-dimensional state spaces is a central challenge in scientific computing, arising prominently in dynamic programming, where the value function of a $d$-dimensional optimal control problem satisfies a Hamilton--Jacobi--Bellman
(HJB) equation over a subset of $\mathbb{R}^d$. Classical methods, such as finite differences \cite{leveque2007finite}, finite elements \cite{zienkiewicz1977finite}, and spectral methods \cite{shen2011spectral}, are effective in low dimensions but typically require meshes, grids, or tensor-product approximation spaces whose cardinality grows exponentially with $d$, a phenomenon known as the curse of dimensionality \cite{bellman1961adaptive}. 
Sparse grids \cite{bungartz2004sparse}, reduced basis methods
\cite{quarteroni2015reduced}, kernel-based approximation
\cite{wendland2004scattered}, and deep learning \cite{han2018solving}
mitigate this growth with ansatze whose complexity is not tied to a grid in $\mathbb{R}^d$. For methods built on a dictionary of candidate features, their effectiveness relies on identifying which basis functions and training samples are most relevant to the solution at hand. 

The prototypical problem motivating this work is the synthesis of optimal feedback laws for nonlinear dynamical systems. Consider a control-affine system with state $\mathbf{y}(t)\in\mathbb{R}^d$ and control signal $\mathbf{u}(t)\in\mathbb{R}^m$, for which we seek a control minimizing an accumulated running cost over an infinite horizon:
\begin{equation}
\begin{aligned}
  &\min _{\mathbf{u}(\cdot)} \int_{0}^{\infty} \left( \ell(\mathbf{y}(t))+\beta \| \mathbf{u}(t)\|_2^2 \right) \, dt, \quad \beta>0,\\
  &\text{subject to }\quad\dot{\mathbf{y}}(t)=\mathbf{f}(\mathbf{y}(t))+\mathbf{q}(\mathbf{y}(t)) \mathbf{u}(t), \quad \mathbf{y}\left(0\right)=\mathbf{x}.
\end{aligned}
\label{eq:infinite_horizon}
\end{equation}
The optimal cost of \eqref{eq:infinite_horizon}, understood as a function of the initial state $\mathbf{x}$, defines the value function $V(\mathbf{x})$ of the control problem. Dynamic programming establishes that, wherever differentiable, the value function satisfies the stationary, first-order Hamilton--Jacobi--Bellman equation
\begin{equation}
0=\min _{\mathbf{u} \in \mathbb{R}^m}\left\{\ell(\mathbf{x})+\beta\|\mathbf{u}\|_2^2+\nabla V(\mathbf{x})^\top \left(\mathbf{f}( \mathbf{x})+\mathbf{q}( \mathbf{x}) \mathbf{u}\right)\right\}. 
\label{eq:stationary}
\end{equation}
Once $V$ is known, the optimal control is recovered in \emph{feedback form}, $\mathbf{u}^*(\mathbf{x})=-\frac{1}{2\beta}\mathbf{q}(\mathbf{x})^\top\nabla V(\mathbf{x})$, valid for every state $\mathbf{x}$. Inserting $\mathbf u^*$ into \eqref{eq:stationary} gives the fully nonlinear PDE $\mathcal L_{\mathrm{HJB}}(V)=0$, with
\begin{equation*}
\mathcal L_{\mathrm{HJB}}(V):=-\frac{1}{4\beta}\nabla V^\top\mathbf q\mathbf q^\top\nabla V+\nabla V^\top\mathbf f+\ell,
\end{equation*}
cast over the state space of the dynamics, so that grid-based discretizations become unfeasible already for moderate values of $d$. The construction of numerical schemes for \eqref{eq:stationary} which remain tractable in high dimensions is the main purpose of this paper. The methodology applies to a general class of first-order PDEs $\mathcal L(V)=0$ (see Section~\ref{sec:model}), with the HJB equation as our leading example.

\paragraph{\textbf{Sparse greedy approximation}}
To approximate the value function $V$, we introduce a large pool of globally defined candidate basis functions $\Phi_{\Theta}(\mathbf{x})=\{\phi_j(\mathbf{x}):j\in\Theta\}$, with $|\Theta|=P\gg 1$, such as multivariate polynomials, and consider approximations of the form  
\begin{equation}
\hat{V}(\mathbf{x};\mathbf{c}) = \sum_{j \in \Gamma} c_j \phi_j(\mathbf{x}) = \mathbf{c}^\top \Phi_{\Gamma}(\mathbf{x}),
\label{eq:approximation}
\end{equation}
where $\Gamma \subset\Theta$ is the index set of active basis functions, $\Phi_{\Gamma}(\mathbf{x}) = [\phi_j(\mathbf{x})]_{j\in\Gamma}^\top$ is the vector of active basis functions, and $\mathbf{c} = [c_j]_{j\in\Gamma}^\top$ is the coefficient vector. In practice only a relatively small subset of $\Phi_{\Theta}$ contributes significantly to the approximation, while the cost of computing $\mathbf c$ grows rapidly with $n:=|\Gamma|$. Our aim is therefore to identify an accurate active set with $n\ll P$. We construct it greedily: starting from a small set $\Gamma^{(1)}$, at iteration $k$ we compute the current approximation $\hat V^{(k)}=\mathbf c^{(k)\top}\Phi_{\Gamma^{(k)}}$, assign to each remaining candidate $\phi_j$, $j\in\Lambda^{(k)}:=\Theta\setminus\Gamma^{(k)}$, an indicator $e_j^{(k)}\ge0$ that predicts the improvement obtained by adding $\phi_j$, and set
\begin{equation*}
\Gamma^{(k+1)}=\Gamma^{(k)}\cup\{j^*\},\qquad j^*\in\operatorname{argmax}_{j\in\Lambda^{(k)}}e_j^{(k)},
\end{equation*}
after which all coefficients are recomputed. The method is thus specified by how $\mathbf c^{(k)}$ is computed and by the choice of $e_j^{(k)}$. We consider two complementary settings.
\begin{itemize}
    \item[i)] \textbf{Model-based framework.} 
    When data are unavailable or expensive to obtain, we compute the coefficients by a Galerkin projection of the PDE and select basis functions with the residual-decay indicator
\begin{equation*}
e_{j,\mathrm{PDE}}^{(k)}=\frac{\bigl|\langle\mathcal L(\hat V^{(k)}),\mathcal L_k\phi_j\rangle_{L^2(\Omega)}\bigr|}{\|\mathcal L_k\phi_j\|_{L^2(\Omega)}},\qquad \mathcal L_k:=D\mathcal L(\hat V^{(k)}).
\end{equation*}
Its square, divided by two, is the largest decrease of the linearized residual functional $\alpha\mapsto\frac12\|\mathcal L(\hat V^{(k)})+\alpha\mathcal L_k\phi_j\|_{L^2(\Omega)}^2$ along $\phi_j$. We derive it from the dual-weighted residual (DWR) framework by projecting the adjoint onto $\operatorname{span}\{\mathcal L_k\phi_j\}$, which makes it closed-form and adjoint-free: it requires only the current residual and the linearized operator, and no mesh.
    \item[ii)] \textbf{Data-driven framework.} Following data-driven approaches to PDE approximation \cite{rudy2019data}, assume we are given a dataset 
$\mathcal D := \{(\mathbf{x}_i, v_i, \mathbf{g}_i)\}_{i=1}^N$ of function values $v_i=V(\mathbf{x}_i)$ and gradients $\mathbf{g}_i=\nabla V(\mathbf{x}_i)$. For HJB equations, such data can be generated in a causality-free fashion by solving open-loop optimal control problems via Pontryagin's Maximum Principle (PMP), see Section~\ref{sec:ocp}. We fit \eqref{eq:approximation} by a gradient-augmented Lasso regression on an active subset of $\mathcal D$. Given prescribed weights $\omega_1,\omega_2>0$ for the value
and gradient errors, at each iteration we add both the sample with the largest misfit
\begin{equation*}
\omega_1|v_i-\hat V^{(k)}(\mathbf x_i)|^2+\omega_2\|\mathbf g_i-\nabla\hat V^{(k)}(\mathbf x_i)\|_2^2
\end{equation*}
and the basis function with the largest supervised-loss indicator $e_{j,\mathrm{SL}}^{(k)}$, which measures the exact one-step decrease of the weighted gradient-augmented least-squares loss along $\phi_j$.   In this way, only a small subset of the available samples is used.
\end{itemize}

\paragraph{\textbf{Related literature}}
Two choices central to this paper are the selection of the approximation space and the states used to construct the approximation. Existing approaches to high-dimensional HJB equations address these choices through different numerical formulations. A first
family determines the approximation from the PDE itself. Polynomial Galerkin schemes
combined with policy iteration \cite{Kalise_2018}, sparse-grid discretizations
\cite{garcke2017suboptimal,bokanowski2013adaptive}, low-rank tensor formats
\cite{dolgov2021tensor,DKS}, and neural networks trained by minimizing the HJB
residual, either pointwise as in physics-informed and Deep Galerkin methods
\cite{Sirignano} or over monotone finite-difference schemes \cite{ETM}, fall into
this class, as do backward SDE formulations \cite{han2018solving}. In these methods
the approximation space, that is, the polynomial degree, the grid level, the tensor
format or the network architecture, is fixed before the solve, or is refined by local
error indicators tied to a grid, as in adaptive sparse grids. 

A second family avoids discretizing the state space by working at individual states
or along optimal trajectories. Since optimal trajectories are characteristics of the
HJB equation, solving open-loop optimal control problems via the PMP yields the value
function and its gradient at individual states without a spatial grid
\cite{kang2015causality,kang2017mitigating,Yegorov,subbotina2006method}, and for
certain Hamiltonians,  Lax--Hopf formulas reduce such pointwise evaluations to convex
optimization problems \cite{DO}. These samples are then used to fit a global
approximation, such as a sparse polynomial \cite{azmi2021optimal} or a neural network
trained with adaptive sampling \cite{nakamura2021adaptive}; since each sample requires
the solution of a two-point boundary value problem, the number of samples largely
determines the cost. Related approaches regress on data from
state-dependent Riccati equations \cite{albi2021gradient}, train networks on the
control objective with HJB penalties  along trajectories \cite{onken2022neural,LVR}, apply symplectic
networks to the PMP system \cite{MZDK}, or learn the feedback law directly
\cite{kunisch2021semiglobal}. 

Adaptive strategies have been developed in several related contexts.
In data-driven settings,
\cite{schaback2019greedy} proposes a greedy point selection algorithm for kernel-based solutions of Dirichlet problems. In machine learning, \cite{siegel2023greedy} develops a greedy training algorithm for shallow neural networks that incrementally adds neurons (basis functions) to minimize empirical risk. In optimal control, \cite{ehring2024hermite,ehring2022greedy}
 combine greedy selection of initial states for open-loop solves with greedy kernel approximation, including Hermite interpolation based on value and gradient data. Their subsequent work \cite{ehring2025online} introduces online enrichment driven by the HJB residual. In these kernel constructions, the approximation space is generated by kernel functions associated with the selected interpolation centers. 
 
In model-based settings, adaptivity is often guided by a posteriori error estimation. A well-established framework is the DWR method, originally developed for finite element methods, which drives adaptive mesh refinement to achieve target accuracy in a prescribed quantity of interest (QoI) \cite{bruchhauser2020dual,becker2001optimal,giles2002adjoint,helml2023plain,oden2001goal}. DWR provides a first-order error approximation for the QoI by weighting the primal residual with the solution of a (linearized) adjoint problem, typically computed with the same numerical scheme as the primal one, such as finite elements, discontinuous Galerkin \cite{cockburn2022adjoint}, finite volumes \cite{chen2014goal}, or feedforward neural networks \cite{roth2022neural}. 
It has been applied to many low-dimensional problems \cite{jayasinghe2018space,venditti2003anisotropic,rabizadeh2020pointwise} and extended to neural networks \cite{fidkowski2021adjoint}, reduced-order models \cite{blonigan2023evaluation}, multiscale \cite{chung2016goal} and virtual element methods \cite{sellmann2025goal}. However, these formulations remain mesh-based, and their effectiveness depends critically on the accuracy of the adjoint approximation \cite{nochetto2009safeguarded}. The approach in \cite{blonigan2023evaluation} projects the adjoint onto a reduced subspace, but still requires solving a reduced adjoint system. Adjoint-based adaptivity for HJB equations has so far been demonstrated only in one- or two-dimensional settings \cite{sagona2003posteriori,cagnetti2013convergence,brenner2021adaptive,gallistl2019mixed}.

 The proposed greedy enrichment is closely related to Orthogonal Matching Pursuit (OMP)~\cite{pati1993orthogonal,Tropp2004} and, more broadly, to greedy dictionary approximation \cite{MallatZhang1993,BinevEtAl2011}. OMP selects the dictionary element with the largest normalized correlation with the current residual and then refits by least squares over the span of the selected elements. The supervised-loss indicator is a weighted OMP selection rule on the gradient-augmented feature matrix; our data-driven algorithm differs from OMP through the greedy selection of training data and the sparsity-promoting refit. The residual-decay indicator is an OMP-type rule in function space, with dictionary $\{\mathcal L_k\phi_j\}$, the PDE residual in place of the data residual, and a Galerkin refit.

\subsection*{Contributions}
The main contributions of this paper are as follows.
\begin{itemize}
  \item[i)] We derive the closed-form residual-decay indicator $e_{j,\mathrm{PDE}}^{(k)}$ for model-based greedy basis enrichment, show that it predicts the first-order reduction of the PDE residual obtained by adding a single candidate basis function, and specialize it to the HJB operator. It is inspired by the goal-oriented DWR framework, but applies to high-dimensional, mesh-free settings while bypassing any adjoint solve. It also admits an empirical version computed from samples.
\item[ii)] We develop a data-driven greedy algorithm for joint basis and training sample
selection, which couples a weighted OMP step for the basis with a greedy selection of gradient-augmented training samples.
\item[iii)] 
We demonstrate the proposed framework on benchmark PDEs and nonlinear optimal control problems in up to $35$ dimensions. The numerical results show that the greedy strategies achieve comparable accuracy to nonadaptive baselines while using substantially fewer basis functions and, in the data-driven setting, fewer
training samples.
\end{itemize}
The rest of the paper is organized as follows. Section~\ref{sec:ocp} presents the finite-horizon optimal control problem, the associated HJB equation, its link with the PMP, and the large-horizon recovery of the infinite-horizon value function. Section~\ref{sec:model} introduces the polynomial dictionaries and the Galerkin and supervised-learning schemes for the coefficients. Sections~\ref{sec:free} and~\ref{sec:data} develop the model-based and data-driven greedy algorithms, and Section~\ref{sec:numerical} presents the numerical experiments.

\section{Optimal control and Hamilton--Jacobi--Bellman PDEs}
\label{sec:ocp}
A primary motivation of this work is the approximation of nonlinear Hamilton-Jacobi-Bellman PDEs arising in optimal control, where the dimensionality make grid-based discretizations intractable. For an initial time $t\in[0,T]$ and an initial state $\mathbf{x}\in\mathbb{R}^d$, we consider a deterministic nonlinear optimal control problem of the form
\begin{equation}
\begin{aligned}
  &\min _{\mathbf{u}(\cdot) \in L^2(t, T ; \mathbb{R}^m)} J\left(\mathbf{u} ; t, \mathbf{x}\right):=\int\limits_{t}^{T}  \ell(\mathbf{y}(s))+\beta \| \mathbf{u}(s)\|_2^2  \; ds+\ell_T(\mathbf{y}(T)),\quad\beta>0,\\
  &\text{subject to }\quad\dot{\mathbf{y}}(s)=\mathbf{f}(\mathbf{y}(s))+\mathbf{q}(\mathbf{y}(s)) \mathbf{u}(s), \quad \mathbf{y}\left(t\right)=\mathbf{x}.
\end{aligned}
\label{eq:optimal_cost}
\end{equation}
Assume that the running cost $\ell: \mathbb{R}^d \rightarrow \mathbb{R}$, the terminal penalty $\ell_T: \mathbb{R}^d \rightarrow \mathbb{R}$, the dynamics $\mathbf{f}: \mathbb{R}^d\rightarrow \mathbb{R}^d$, and the control-to-state-operator $\mathbf{q}: \mathbb{R}^d \rightarrow \mathbb{R}^{d \times m}$ are continuously differentiable  in $\mathbb{R}^d$. The value function of \eqref{eq:optimal_cost} is defined as the optimal cost-to-go
\begin{equation*}
V(t,\mathbf{x}):=\min _{\mathbf{u}(\cdot) \in L^2(t, T ; \mathbb{R}^m)} J\left(\mathbf{u} ; t, \mathbf{x}\right),\qquad (t,\mathbf{x})\in[0,T]\times\mathbb{R}^d.
\end{equation*}
By the Dynamic Programming Principle \cite{bertsekas2012dynamic}, wherever $V$ is differentiable it satisfies the time-dependent HJB equation
\begin{equation*}
\left\{\begin{array}{l}
\partial_t V(t,\mathbf{x})+\min\limits_{\mathbf{u}\in\mathbb{R}^m}\left\{\ell(\mathbf{x})+\beta\|\mathbf{u}\|_2^2+\nabla V(t,\mathbf{x})^\top\left(\mathbf{f}(\mathbf{x})+\mathbf{q}(\mathbf{x})\mathbf{u}\right)\right\}=0,\\
V(T,\mathbf{x})=\ell_T(\mathbf{x}).
\end{array}\right.
\end{equation*}
from where the optimal control in feedback form is given by
\begin{equation*}
\mathbf{u}^*(t, \mathbf{x})=-\frac{1}{2 \beta} \mathbf{q}(\mathbf{x})^\top \nabla V(t, \mathbf{x}),
\end{equation*}
yielding the equivalent HJB PDE
\begin{equation}
\label{finite_HJB}
\left\{\begin{aligned}
&\partial_t V(t,\mathbf{x})+\nabla V(t,\mathbf{x})^\top \mathbf{f}(\mathbf{x})-\frac{1}{4 \beta}\left\|\mathbf{q}(\mathbf{x})^\top \nabla V(t,\mathbf{x})\right\|_2^2+\ell(\mathbf{x})=0, \\
&V(T, \mathbf{x})=\ell_T(\mathbf{x}).
\end{aligned}\right.
\end{equation}
An alternative solution to the optimal control problem \eqref{eq:optimal_cost} follows Pontryagin's Maximum Principle, corresponding to first-order necessary optimality conditions in the form of a two-point boundary-value problem (TPBVP) coupling the optimal state, adjoint, and control $\left(\mathbf{y}^*(s), \mathbf{p}^*(s), \mathbf{u}^*(s)\right)$ over $s\in(t,T)$:
\begin{equation*}
\left\{\begin{aligned}
\frac{d}{ds} \mathbf{y}^*(s) & =\mathbf{f}\left(\mathbf{y}^*(s)\right)+\mathbf{q}\left(\mathbf{y}^*(s)\right) \mathbf{u}^*(s), \\
\mathbf{y}^*\left(t\right) & =\mathbf{x}, \\
-\frac{d}{ds} \mathbf{p}^*(s) & =\partial_{\mathbf{y}}\left(\mathbf{f}\left(\mathbf{y}^*(s)\right)+\mathbf{q}\left(\mathbf{y}^*(s)\right) \mathbf{u}^*(s)\right)^\top \mathbf{p}^*(s)+\partial_{\mathbf{y}} \ell\left(\mathbf{y}^*(s)\right), \\
\mathbf{p}^*(T) & =\nabla \ell_T\left(\mathbf{y}^*(T)\right),
\end{aligned}\right.
\end{equation*}
closed via the optimal control
\begin{equation*}
\mathbf{u}^*(s)=-\frac{1}{2 \beta} \mathbf{q}\left(\mathbf{y}^*(s)\right)^\top \mathbf{p}^*(s), \quad \forall s \in\left(t, T\right).
\end{equation*}
Under the regularity assumption $V\in C^1$, the PMP forward--backward dynamics act as characteristic curves of the HJB equation \cite{subbotina2006method,zhou1990maximum,pontryagin2018mathematical}. Thus,  the value function at $(t,\mathbf{x})$ can be computed 
by solving the TPBVP and evaluating the optimal cost \eqref{eq:optimal_cost} along the resulting optimal trajectory, while its gradient is identified as the adjoint, $\nabla V(t,\mathbf{x})=\mathbf{p}^*(t)$.
Consequently, the data-driven approach developed in \cite{azmi2021optimal}, which constitutes one of the building blocks of the present work, approximates the solution of \eqref{finite_HJB} by first generating a gradient-augmented dataset using a characteristic-based, causality-free method derived from solving the TPBVP via PMP, and then applies sparse regression to recover the value function using a polynomial basis. Building upon this framework, we introduce greedy selection strategies for both basis functions and data points. Instead of prescribing a fixed basis and dataset, we iteratively enrich the approximation subspace by selecting the most informative basis functions and sampling data points. This adaptive strategy provides accurate approximations of value functions while improving computational efficiency.

\noindent We conclude this section by connecting the finite-horizon problem \eqref{eq:optimal_cost} to the infinite-horizon problem \eqref{eq:infinite_horizon} and its stationary HJB equation \eqref{eq:stationary}. Note that \eqref{eq:infinite_horizon} cannot be treated as the instance $T=\infty$ of \eqref{eq:optimal_cost}: without additional admissibility assumptions on the class of controls, the resulting problem is ill-posed. Instead, we exploit the fact that, for stabilizable dynamics and a coercive running cost, the finite-horizon value function becomes approximately stationary sufficiently far away from the terminal time, and the influence of the terminal penalty $\ell_T$ on $V(0,\cdot)$ fades as $T$ grows. Hence, for $T$ sufficiently large, the stationary value function $V(\mathbf{x})$ and the associated feedback law $\mathbf{u}^*(\mathbf{x})$ are recovered from their time-dependent counterparts through
\begin{equation*}
V(\mathbf{x})\approx V(0,\mathbf{x}),\qquad \mathbf{u}^*(\mathbf{x})\approx \mathbf{u}^*(0,\mathbf{x}).
\end{equation*}
This large-horizon approximation is the route followed throughout the paper: in the model-based setting, the time-dependent HJB equation \eqref{finite_HJB} is integrated backwards from the terminal condition and $V(0,\cdot)$ is retained, while in the data-driven setting, PMP-based samples of $V(0,\cdot)$ and $\nabla V(0,\cdot)$ are used as training data.

\section{Polynomial models and approximation of HJB PDEs}
\label{sec:model}
The adaptive strategies developed in this paper apply to general stationary problems of the form
\begin{equation}
(\mathcal{L}V)(\mathbf{x}) := \mathcal{H}\big(\mathbf{x},\,V(\mathbf{x}),\,\nabla V(\mathbf{x})\big)=0,
\qquad \mathbf{x}\in\Omega,
\label{eq:pde}
\end{equation}
where $\Omega\subset\mathbb{R}^d$ is a compact domain, $\mathcal{V}$ is a Hilbert
space continuously embedded in $H^1(\Omega)$, and $\mathcal{L}:\mathcal{V}\to
L^2(\Omega)$ is the (possibly nonlinear) first-order differential operator induced
by the Hamiltonian $\mathcal{H}:\Omega\times\mathbb{R}\times\mathbb{R}^d\to\mathbb{R}$.
Our leading example is the stationary HJB equation \eqref{eq:stationary}, for which $\mathcal{H}$ is the control Hamiltonian obtained upon minimization over $\mathbf{u}$, and $\mathcal{L}$ is the stationary HJB operator made precise in Section~\ref{sec:HJB_indicator}. Recall that we use $V$ to denote the solution of \eqref{eq:pde}, and denote its approximation by $\hat{V}$.

\noindent In this section, we introduce the construction of the candidate basis set $\Phi_\Theta$ and present two complementary numerical schemes for the approximation of the PDE \eqref{eq:pde}. In the model-based setting, the coefficients are computed by a Galerkin projection of the PDE operator $\mathcal{L}$. In the data-driven setting, the coefficients are obtained by supervised learning from gradient-augmented samples of the value function.

\subsection{Sparse polynomial ansatz}
Let $\Omega=\prod_{i=1}^d\Omega_i$ be a bounded domain, where $\Omega_i\subset \mathbb{R}$. For $i=1,\cdots,d$, let $\{\phi^{(i)}_{m}\}_{m\in\mathbb{N}_0}$ be a one-dimensional orthonormal basis
of $L^2(\Omega_i)$. Typical choices include orthonormal polynomial basis, such
as Legendre and orthonormal trigonometric basis. % such as normalized sine and cosine functions. % smooth multiresolution families such as  Symlets or Meyer wavelet \cite{mallat1999wavelet}.
Then for $j=(j_1,j_2,\ldots,j_d)\in\mathbb{N}_0^d$, we define the corresponding tensor-product basis function by 
\begin{equation*}
\phi_j(\mathbf{x})=\prod_{i=1}^d \phi^{(i)}_{j_i}\left(x_i\right),\qquad j=(j_1,j_2,\ldots,j_d)\in\mathbb{N}_0^d.
\end{equation*}
To manage complexity in high dimensions, we use the total degree index set
$$
\text{(Total degree)}\quad \Theta^{\mathrm{TD}}(M):=\left\{j\in\mathbb{N}_0^d: \;\sum_{i=1}^d j_i \leq M\right\},$$ and the hyperbolic cross index set
\begin{align*}
    \text{(Hyperbolic cross)}\quad  \Theta^{\mathrm{HC}}(\mathfrak{J}):=\left\{j\in\mathbb{N}_0^d: \;\prod_{i=1}^d\left(j_i+1\right) \leq 2^\mathfrak{J}+1\right\},
\end{align*}
where $M,\mathfrak{J}\in\mathbb{N}$. For comparison, the full tensor index set  $$  \text{(Full tensor)}\quad \Theta^{\mathrm{FULL}}(S):=\left\{j\in\mathbb{N}_0^d: \;\max_{i=1,\ldots,d} j_i \leq S\right\}$$
with $S\in\mathbb{N}$ has cardinality
$|\Theta^{\mathrm{FULL}}(S)|=(S+1)^d$, which grows
exponentially with $d$. By contrast,  the total-degree index set has cardinality
$|\Theta^{\mathrm{TD}}(M)|=\sum_{m=0}^M\binom{d+m-1}{m}$, which grows combinatorially with $M$ and $d$. The hyperbolic cross set provides a further reduction and thus scales better for high-dimensional problems, with  cardinality satisfying $|\Theta^{\mathrm{HC}}(\mathfrak{J})|\leq \min \left\{2^{3\mathfrak{J}+1}4^d,\; e^2 (4d)^\mathfrak{J}\right\}$ \cite{azmi2021optimal}.
%\begin{remark}
%As an alternative to the separable tensor-product constructions considered above, we may construct the candidate set $\Phi$ using localized Gaussian (radial basis) functions. Such a basis is easy to generate in arbitrary dimension $d$ (by choosing centers and widths), but they are not orthogonal in $L^2(\Omega)$, which can lead to ill-conditioned mass matrices.
%\end{remark}
\subsection{Model-based schemes: Galerkin projection}
Given a selected basis index set $\Gamma \subset \Theta$ with
$|\Gamma| = n$, we define the  $n$-dimensional approximation subspace
\begin{equation*}
\mathcal{V}_{n}:=\operatorname{span}\left\{\phi_j: j \in \Gamma\right\} \subset \mathcal{V}.
\end{equation*}
Writing $\hat{V}=\Phi_\Gamma^\top\mathbf{c}$ with $\mathbf{c}=[c_j]_{j\in\Gamma}\in\mathbb{R}^n$, the Galerkin approximation of \eqref{eq:pde} in $\mathcal{V}_n$ is determined by the $n$ equations
\begin{equation}
\mathbf{H}(\mathbf{c})=0,\qquad \mathbf{H}_i(\mathbf{c}):=\left\langle\mathcal{L}\big(\Phi_\Gamma^\top\mathbf{c}\big), \phi_i\right\rangle_{L^2(\Omega)},\quad i\in\Gamma.
\label{eq:galerkin}
\end{equation}
When $\mathcal{L}$ is affine, \eqref{eq:galerkin} is a linear system for $\mathbf{c}$; this is the case of the benchmark in Section~\ref{sec:testcos}. The ansatz and the projection \eqref{eq:galerkin} are formulated for a stationary PDE. The HJB equation \eqref{finite_HJB} associated with the control problem of Section~\ref{sec:ocp}, however, is time-dependent: we approximate its solution $V(t,\cdot)$ on $[0,T]$ and retain $V(0,\cdot)$, which approximates the stationary value function for large $T$. We therefore extend the ansatz to time-dependent functions by keeping the spatial basis fixed and letting the coefficients depend on time. Projecting in space as in \eqref{eq:galerkin} then yields a system of ordinary differential equations for the coefficients, which we now derive.
\subsubsection{A time-dependent Galerkin scheme for HJB equations}
\label{sec: finite}
The value function of the finite-horizon problem \eqref{eq:optimal_cost} satisfies \eqref{finite_HJB}, which in terms of the stationary HJB operator $\mathcal{L}_{\mathrm{HJB}}$ of \eqref{eq:HJBoperator} reads $\partial_t V+\mathcal{L}_{\mathrm{HJB}}V=0$ with $V(T,\cdot)=\ell_T$. As in \eqref{eq:approximation}, we approximate the value function by
\begin{equation*}
\hat{V}(t,\mathbf{x}) =\sum_{j\in\Gamma} c_j(t) \phi_j(\mathbf{x})=\Phi_\Gamma(\mathbf{x})^\top\mathbf{c}(t),
\end{equation*}
where $\mathbf{c}(t) = [c_j(t)]_{j\in\Gamma}^\top$ is the vector of time-dependent coefficients. Projecting the HJB equation onto the span of $\{\phi_i\}_{i \in \Gamma}$, for each $i\in\Gamma$, we solve for $\mathbf{c}(t)$ such that
\begin{equation*}
\left\langle\partial_t V(t,\mathbf{x})+\nabla V(t,\mathbf{x})^\top \mathbf{f}(\mathbf{x})-\frac{1}{4 \beta}\left\|\mathbf{q}(\mathbf{x})^\top \nabla V(t,\mathbf{x})\right\|_2^2+\ell(\mathbf{x}), \phi_i(\mathbf{x})\right\rangle_{L^2(\Omega)}=0.
\end{equation*}
Substituting $\hat{V}$  into the above equation,
we obtain
%\begin{align*}\left\langle \sum_{j\in\Gamma} \dot{c}_j(t) \phi_j, \phi_i \right\rangle +\left\langle \sum_{j\in\Gamma}c_j(t) \nabla \phi_j^{\top} \mathbf{f}, \phi_i \right\rangle- \frac{1}{4\beta} \left\langle \left| \mathbf{q}^\top \left(\sum_{j\in\Gamma} c_j(t) \nabla \phi_j\right)\right|^2, \phi_i \right\rangle + \left\langle \ell, \phi_i\right\rangle = 0\end{align*}
\begin{equation}
\begin{aligned}
&\left\langle \sum_{j\in\Gamma} \dot{c}_j(t) \phi_j(\mathbf{x}), \phi_i(\mathbf{x}) \right\rangle_{L^2(\Omega)} +
\left\langle \sum_{j\in\Gamma}c_j(t) \nabla \phi_j(\mathbf{x})^\top \mathbf{f}(\mathbf{x}), \phi_i(\mathbf{x}) \right\rangle_{L^2(\Omega)}\\ 
&- \frac{1}{4\beta} \left\langle \left\| \mathbf{q}(\mathbf{x})^\top \left(\sum_{j\in\Gamma} c_j(t) \nabla \phi_j(\mathbf{x}) \right)\right\|_2^2, \phi_i(\mathbf{x}) \right\rangle_{L^2(\Omega)} 
+ \left\langle \ell(\mathbf{x}), \phi_i(\mathbf{x}) \right\rangle_{L^2(\Omega)} = 0.
\end{aligned}
\label{eq:galerkin_hjb_discrete}
\end{equation}
This motivates the definition of the following matrices and vectors:
\begin{align*}
\mathbf{M}\in\mathbb{R}^{n\times n}\quad &\text{with}\quad\mathbf{M}_{i j}=\int_\Omega \phi_i(\mathbf{x}) \phi_j(\mathbf{x}) d \mathbf{x},\\   
\mathbf{K}\in\mathbb{R}^{n\times n}\quad &\text{with}\quad \mathbf{K}_{ij}=\int_{\Omega} \phi_i(\mathbf{x}) \nabla \phi_j(\mathbf{x})^\top \mathbf{f}(\mathbf{x}) d \mathbf{x},\\
\mathbf{L}\in\mathbb{R}^{n}\quad &\text{with}\quad \mathbf{L}_{i}=\int_{\Omega} \phi_i(\mathbf{x}) \ell(\mathbf{x}) d \mathbf{x},\\
\mathbf{L}_T\in\mathbb{R}^{n}\quad &\text{with}\quad \left(\mathbf{L}_T\right)_{i}=\int_{\Omega} \phi_i(\mathbf{x}) \ell_T(\mathbf{x}) d \mathbf{x},\\
\mathbf{Q}(\mathbf{c})\in\mathbb{R}^{n}\quad &\text{with}\quad \mathbf{Q}_i(\mathbf{c})=\int_{\Omega}\phi_i(\mathbf{x})\left\|\mathbf{q}(\mathbf{x})^\top\left(\sum_{j \in \Gamma} c_j \nabla \phi_j(\mathbf{x})\right)\right\|_2^2 d \mathbf{x}.
\end{align*}
Since
$\partial_t \hat{V}(t, \mathbf{x})=\sum_{j\in\Gamma} \dot{c}_j(t) \phi_j(\mathbf{x}),$
we have 
\begin{equation*}
\left\langle\partial_t\hat{V}, \phi_i\right\rangle_{L^2(\Omega)}=\sum_{j\in\Gamma}\dot{c}_j(t)\left\langle\phi_j, \phi_i\right\rangle_{L^2(\Omega)}=[\mathbf{M} \dot{\mathbf{c}}(t)]_i.
\end{equation*}
Finally, we obtain the system
\begin{equation}\dot{\mathbf{c}}(t)=\mathbf{M}^{-1}\left(-\mathbf{K}\mathbf{c}(t)+\frac{1}{4 \beta} \mathbf{Q}(\mathbf{c}(t))-\mathbf{L}\right).
\label{eq:ode_t}\end{equation}
We solve \eqref{eq:ode_t} backward in time on $[0,T]$ with terminal condition $\mathbf{M}\mathbf{c}(T)=\mathbf{L}_T$, i.e., $\hat V(T,\cdot)$ is the $L^2(\Omega)$ projection of $\ell_T$ onto $\mathcal{V}_n$ (in particular, $\mathbf{c}(T)=0$ when $\ell_T\equiv 0$). Let $F(\mathbf{c}(t))$ denote the right-hand side of \eqref{eq:ode_t}. %We consider two options for solving the ODE system \eqref{eq:ode_t} backward in time.
In an explicit Euler scheme for backward integration, we discretize $[0, T]$ into $N_t$ time steps of size $\Delta t=T/{N_t}$, and define the reversed time $t_{n_t}=T-n_t \;\Delta t$, where $n_t=0, \cdots, N_t$. Then, starting from $\mathbf{c}_0=\mathbf{c}(T)$, we update  $$\mathbf{c}_{n_t+1}=\mathbf{c}_{n_t}-\Delta t \;F(\mathbf{c}_{n_t}),\quad n_t= 0, \cdots, N_t-1,$$
and the value function at time
$t=0$ is approximated as
\begin{equation*}
\hat V(0,\mathbf{x}) \approx\sum_{j\in\Gamma} c_j(0) \phi_j(\mathbf{x}).
\end{equation*}
%This approximation also relates to the nonlinear model predictive control, where an open-loop optimal control problem is solved over a finite horizon, and only the initial control input $\mathbf{u}^*(0)$ is applied to evolve the state for a short period. The process is then repeated from the updated state, reflecting the idea of using initial-time value function or control information to guide system evolution.

\subsection{Data-driven schemes: supervised learning with gradient-augmented data}
\label{sec:sp}
We consider a training dataset $\mathcal{D}=\{(\mathbf{x}_i, v_i, \mathbf{g}_i)\}_{i=1}^N$ and an active subset of the basis $\Phi_\Gamma=\{\phi_j\}_{j\in\Gamma}$, with $|\Gamma|=n$.
For each sample $\mathbf{x}_i$,  we define the feature row corresponding to
function values as  $\Phi_\Gamma(\mathbf{x}_i):=\{\phi_j(\mathbf{x}_i)\}_{j\in\Gamma}\in\mathbb{R}^n$. For the gradient information, we define, for each coordinate
$k = 1,\dots,d$,
\begin{equation*}
\Phi_\Gamma^{(k)}\left(\mathbf{x}_i\right):=\left\{
\frac{\partial \phi_j\left(\mathbf{x}_i\right)}{\partial x_k}
\right\}_{j\in\Gamma}\in\mathbb{R}^{n}
\end{equation*}
and the full gradient feature matrix as
\begin{equation*}
\mathbf{G}(\mathbf{x}_i)=\left[\begin{array}{c}
\Phi_\Gamma^{(1)}\left(\mathbf{x}_i\right) \\
\Phi_\Gamma^{(2)}\left(\mathbf{x}_i\right) \\
\vdots \\
\Phi_\Gamma^{(d)}\left(\mathbf{x}_i\right)
\end{array}\right]\in\mathbb{R}^{ d\times n}.
\end{equation*}
The complete feature matrix is then defined by
\begin{equation*}
\mathbf{X}=\left[\begin{array}{c}
\tilde{\mathbf{F}} \\ \tilde{\mathbf{G}}
\end{array}\right]\in\mathbb{R}^{N(1+d)\times n},
\end{equation*}
where 
\begin{equation*}
\tilde{\mathbf{F}} =\left[\begin{array}{c}
\Phi_\Gamma\left(\mathbf{x}_1\right) \\
\Phi_\Gamma\left(\mathbf{x}_2\right) \\
\vdots \\
\Phi_\Gamma\left(\mathbf{x}_{N}\right)
\end{array}\right]\in\mathbb{R}^{N\times n},\quad \tilde{\mathbf{G}}=\left[\begin{array}{c}
\mathbf{G}\left(\mathbf{x}_1\right) \\
\mathbf{G}\left(\mathbf{x}_2\right) \\
\vdots \\
\mathbf{G}\left(\mathbf{x}_{N}\right)
\end{array}\right]\in\mathbb{R}^{N  d\times n}.
\end{equation*}
Similarly,  the corresponding target vector $\mathbf{b}\in\mathbb{R}^{ N\left(1+d\right)}$ is constructed by concatenating the function values and gradients: \begin{equation*}
\mathbf{b}=\left[\begin{array}{c}
\mathbf{v}\\ 
\mathbf{g}
\end{array}\right],\quad\mathbf{v}=\left(v_1,v_2,\cdots,v_N\right)^\top\in\mathbb{R}^{N},\quad\mathbf{g}=\left(\mathbf{g}_1^\top, \cdots, \mathbf{g}_N^\top\right)^\top\in\mathbb{R}^{dN}.
\end{equation*}
The optimal coefficient vector $\mathbf{c}\in\mathbb{R}^n$ is obtained by solving the following LASSO regression problem
\begin{equation*}
\mathbf{c}^*=\operatorname{argmin}_{\mathbf{c} \in \mathbb{R}^n}\left\{\|\mathbf{b}-\mathbf{X}\mathbf{c}\|_2^2+\lambda\|\mathbf{c}\|_1\right\},
\end{equation*}
where $\lambda>0$ is a regularization parameter that controls the sparsity in the solution.

\section{Model-based adaptive basis selection}
\label{sec:free}
In this section, we develop a closed-form indicator that predicts the reduction in the PDE residual obtained by adding a single candidate basis function to the current approximation space. The indicator is inspired by the goal-oriented DWR framework for a suitable choice of QoI, while extending this idea to high-dimensional, mesh-free settings and bypassing any adjoint solve. 

\noindent Recall that our goal is to solve \eqref{eq:pde} with an operator $\mathcal{L}: \mathcal{V} \rightarrow L^2(\Omega)$,
 where  $\mathcal{V}$ is a Hilbert space continuously embedded in $H^1(\Omega)$. 
For any $g\in L^2(\Omega)$ and $\phi\in\mathcal{V}$, the embedding $\mathcal{V} \hookrightarrow L^2(\Omega)\hookrightarrow (\mathcal{V})^*$ allows us to identify the duality pairing $\langle\cdot, \cdot\rangle_{\mathcal{V}^*, \mathcal{V}}$ between  $\mathcal{V}^*$ and  $\mathcal{V}$ with the $L^2$ inner product
\begin{equation*}
\langle g, \phi\rangle_{\mathcal{V}^*, \mathcal{V}}:=\int_{\Omega} g(\mathbf{x}) \phi(\mathbf{x}) d \mathbf{x}=\langle g, \phi\rangle_{L^2(\Omega)}.
\end{equation*}
At the $k$-th iteration, let $\hat V^{(k)}$ denote the
current approximation. We assume that \eqref{eq:pde} admits a solution $V\in\mathcal V$, which is fixed throughout as the
target of the approximation, and that,
for each iterate $\hat V^{(k)}$, the operator $\mathcal L$ is Fréchet
differentiable in an open neighborhood of $\hat V^{(k)}$. This
differentiability is used only locally in the construction of the residual-based
indicator below.
\paragraph{Index-set partition} At the $k$-th iteration, the prescribed basis pool $\Phi_{\Theta}$ is partitioned as \begin{equation*}
\Phi_{\Theta}(\mathbf{x})=\Phi_{\Gamma^{(k)}}(\mathbf{x}) \cup \Phi_{\Lambda^{(k)}}(\mathbf{x}),\quad%\Gamma^{(k)} \cup \Lambda^{(k)}=\Theta,\quad
\Gamma^{(k)} \cap \Lambda^{(k)}=\emptyset,
\end{equation*}
where $\Gamma^{(k)} \subset\Theta$ denotes the index set of selected basis functions, $n_k:=\left|\Gamma^{(k)}\right|$, and $\Lambda^{(k)}=\Theta \backslash \Gamma^{(k)}$ is the index set of remaining candidates. We construct $\hat V^{(k+1)}$ by greedily selecting a single basis function from the remaining candidates $\Lambda^{(k)}$ and adding it to the active set $\Gamma^{(k)}$.

\noindent Given $\Gamma^{(k)}$, let $\mathcal{V}_{n_k}:=\operatorname{span}\left\{\phi_j: j \in \Gamma^{(k)}\right\} \subset \mathcal{V}$ and let $\hat{V}^{(k)} \in \mathcal{V}_{n_k}$ be the current approximation. We define the PDE residual
\begin{equation}
R\left(\hat{V}^{(k)}\right):=\mathcal{L}\left(\hat{V}^{(k)}\right)\in L^2(\Omega).
\label{eq:pde_residual}
\end{equation}
Since $\mathcal{L}$ is Fr\'echet differentiable at $\hat{V}^{(k)}$, its derivative $\mathcal{L}_k:=D \mathcal{L}\left(\hat{V}^{(k)}\right):\mathcal{V} \rightarrow L^2(\Omega)$ is a bounded linear operator. In particular, for all  $\psi \in \mathcal{V}$, the directional derivative
\begin{equation*}
D\mathcal{L}\left(\hat{V}^{(k)}\right)[\psi]=\lim _{\epsilon \rightarrow 0} \frac{\mathcal{L}\left(\hat{V}^{(k)}+\epsilon \psi\right)-\mathcal{L}\left(\hat{V}^{(k)}\right)}{\epsilon}
\end{equation*}
exists and satisfies $D \mathcal{L}\left(\hat{V}^{(k)}\right)[ \psi]=\mathcal{L}_k \psi$. Its adjoint $\mathcal{L}_k^*:L^2(\Omega)\to\mathcal{V}^*$ is characterized by
\begin{equation*}
\langle \mathcal{L}_k^* w,\,v\rangle_{\mathcal{V}^*,\mathcal{V}}=\langle w,\,\mathcal{L}_k v\rangle_{L^2(\Omega)},\quad \forall v\in\mathcal{V},\;w\in L^2(\Omega).
\end{equation*}

\subsection{Standard dual-weighted residual method }
We first briefly recall the standard dual-weighted residual (DWR) method \cite{becker2001optimal}, in which the error in a target quantity of interest (QoI) is approximated by weighting the primal residual with the solution of a linearized adjoint equation \cite[Section~3.1]{blonigan2023evaluation}. This provides the conceptual starting point for our indicator.  In the following, we replace
the expensive adjoint solve with a computable projection tailored to adaptive
basis selection.

\noindent We reformulate \eqref{eq:pde} as a PDE-constrained optimization problem
\begin{equation}
\min _{V \in \mathcal{V}} \mathcal{Q}(V) \quad \text { subject to } \mathcal{L}(V)=0,
\label{eq:constrained_opt}
\end{equation}
where $\mathcal{Q}: \mathcal{V} \rightarrow \mathbb{R}$ is a Fr\'echet differentiable functional that specifies the QoI but does not affect the PDE solution.
The associated Lagrangian is defined by
\begin{equation*}
\mathfrak{L}(V, p):=\mathcal{Q}(V)+\langle p, \mathcal{L}(V)\rangle_{L^2(\Omega)},
\end{equation*}
with adjoint variable $p \in L^2(\Omega)$. First-order optimality  conditions yield the adjoint equation
\begin{equation}
\mathcal{Q}^{\prime}(V)+D\mathcal{L}(V)^*[p]=0\quad\text{ in } \mathcal{V}^*.
\label{eq:adjoint}
\end{equation}
At the current iterate we use the adjoint equation linearized at $\hat V^{(k)}$: let $p_k\in L^2(\Omega)$ satisfy
\begin{equation}
\mathcal{L}_k^*p_k=-\mathcal{Q}'\bigl(\hat V^{(k)}\bigr)\quad\text{ in } \mathcal{V}^*.
\label{eq:adjoint_k}
\end{equation}
We denote by
\begin{align*}
r_{\mathcal L}(h)&:=\mathcal L\bigl(\hat V^{(k)}+h\bigr)-R\bigl(\hat V^{(k)}\bigr)-\mathcal L_k h,\\
r_{\mathcal Q}(h)&:=\mathcal Q\bigl(\hat V^{(k)}+h\bigr)-\mathcal Q\bigl(\hat V^{(k)}\bigr)-\bigl\langle \mathcal Q'\bigl(\hat V^{(k)}\bigr),h\bigr\rangle_{\mathcal V^*,\mathcal V},
\end{align*}
the first-order Taylor remainders of $\mathcal L$ and $\mathcal Q$ at $\hat V^{(k)}$, so that $\|r_{\mathcal L}(h)\|_{L^2(\Omega)}=o(\|h\|_{\mathcal V})$ and $r_{\mathcal Q}(h)=o(\|h\|_{\mathcal V})$ as $h\to0$ in $\mathcal V$. The following lemma is the error representation underlying the DWR method.

\begin{lemma}[DWR error representation]
\label{lem:dwr}
Let $V\in\mathcal V$ solve \eqref{eq:pde}, let $p_k\in L^2(\Omega)$ satisfy \eqref{eq:adjoint_k}, and set $\delta V:=V-\hat V^{(k)}$. Then
\begin{equation}
\begin{aligned}
\mathcal{Q}(V)-\mathcal{Q}\bigl(\hat{V}^{(k)}\bigr) ={}&\bigl\langle R\bigl(\hat{V}^{(k)}\bigr), p_k\bigr\rangle_{L^2(\Omega)}\\
&+r_{\mathcal Q}(\delta V)+\bigl\langle p_k, r_{\mathcal L}(\delta V)\bigr\rangle_{L^2(\Omega)}.
\end{aligned}
\label{eq:proof}
\end{equation}
In particular, the last two terms are $o(\|\delta V\|_{\mathcal V})$, and they vanish if $\mathcal L$ is affine and $\mathcal Q$ is linear.
\end{lemma}
\begin{proof}
Since $\mathcal L(V)=0$, the definition of $r_{\mathcal L}$ with $h=\delta V$ gives $\mathcal L_k\delta V=-R(\hat V^{(k)})-r_{\mathcal L}(\delta V)$. Using the definition of $r_{\mathcal Q}$, the adjoint equation \eqref{eq:adjoint_k} and the definition of $\mathcal L_k^*$, we obtain
\begin{align*}
\mathcal{Q}(V)-\mathcal{Q}\bigl(\hat V^{(k)}\bigr)
&=\bigl\langle \mathcal{Q}'\bigl(\hat V^{(k)}\bigr),\delta V\bigr\rangle_{\mathcal V^*,\mathcal V}+r_{\mathcal Q}(\delta V)\\
&=-\bigl\langle p_k,\mathcal L_k\delta V\bigr\rangle_{L^2(\Omega)}+r_{\mathcal Q}(\delta V)\\
&=\bigl\langle p_k,R\bigl(\hat V^{(k)}\bigr)\bigr\rangle_{L^2(\Omega)}+\bigl\langle p_k,r_{\mathcal L}(\delta V)\bigr\rangle_{L^2(\Omega)}+r_{\mathcal Q}(\delta V).
\end{align*}
The final assertion follows from the bounds on $r_{\mathcal L}$ and $r_{\mathcal Q}$ and the Cauchy--Schwarz inequality. \qed
\end{proof}
The leading term $\langle R(\hat V^{(k)}),p_k\rangle_{L^2(\Omega)}$ in \eqref{eq:proof} is the standard goal-oriented DWR error estimate \cite{becker2001optimal,suttmeier2008numerical}.
The adjoint $ p_k$ serves as a sensitivity kernel, quantifying the influence of the residual $ R(\hat{V}^{(k)})$ on the functional error $\mathcal{Q}(V) - \mathcal{Q}(\hat{V}^{(k)})$.

\noindent In the following, we introduce a computable indicator by projecting the adjoint $p_k$ onto the linearized image of candidate directions.

\subsection{Single-direction projected DWR indicator}
\label{sec:single}
Given $\Gamma^{(k)}$ and $\hat{V}^{(k)}\in \mathcal{V}_{n_k}$, define $\mathcal{V}_{\Lambda^{(k)}}:=\operatorname{span}\left\{\phi_j: j \in \Lambda^{(k)}\right\} \subset \mathcal{V}$ and its linearized image
\begin{equation*}
\mathcal{W}_{\Lambda^{(k)}}:=\operatorname{span}\left\{w_j:=\mathcal{L}_k \phi_j: j \in \Lambda^{(k)}\right\} \subset L^2(\Omega).
\end{equation*}
We project $p_k$ onto the linearized image $\operatorname{span}\{w_j\}\subset\mathcal W_{\Lambda^{(k)}}$ of a single candidate direction $\phi_j$. The next lemma shows that this projection can be computed without solving the adjoint equation.

\begin{lemma}[Projected adjoint]
\label{lem:proj}
Let $j\in\Lambda^{(k)}$ with $w_j=\mathcal L_k\phi_j\neq0$, and let $p_k\in L^2(\Omega)$ satisfy \eqref{eq:adjoint_k}. Then the $L^2(\Omega)$-orthogonal projection of $p_k$ onto $\operatorname{span}\{w_j\}$ is
\begin{equation}
p_k^{\mathrm{proj},j}=\frac{\bigl\langle-\mathcal{Q}^{\prime}\bigl(\hat{V}^{(k)}\bigr), \phi_j\bigr\rangle_{\mathcal{V}^*, \mathcal{V}}}{\left\|w_j\right\|_{L^2(\Omega)}^2} \; w_j.
\label{eq:p_proj}
\end{equation}
\end{lemma}
\begin{proof}
In this proof, $\langle\cdot,\cdot\rangle$ and $\|\cdot\|$ denote the inner product and norm of $L^2(\Omega)$. The orthogonal projection onto $\operatorname{span}\{w_j\}$ maps $p$ to $\langle p,w_j\rangle\|w_j\|^{-2}\,w_j$. By the definition of $\mathcal L_k^*$ and \eqref{eq:adjoint_k},
\begin{equation*}
\langle p_k,w_j\rangle=\langle p_k,\mathcal L_k\phi_j\rangle=\langle \mathcal L_k^*p_k,\phi_j\rangle_{\mathcal V^*,\mathcal V}=\bigl\langle -\mathcal Q'\bigl(\hat V^{(k)}\bigr),\phi_j\bigr\rangle_{\mathcal V^*,\mathcal V}. \qquad\qed
\end{equation*}
\end{proof}
Replacing $p_k$ with $p_k^{\mathrm{proj},j}$ in the leading term of \eqref{eq:proof} gives the single-direction projected DWR indicator
\begin{equation}
\mathcal{E}^{\mathrm{proj}}_j
:= \big\langle R(\hat V^{(k)}),\, p_k^{\mathrm{proj},j}\big\rangle_{L^2(\Omega)}
= \frac{\big\langle -\mathcal{Q}'(\hat V^{(k)}),\,\phi_j\big\rangle_{\mathcal V^*,\mathcal V}
       \big\langle R(\hat V^{(k)}),\, w_j\big\rangle_{L^2(\Omega)}}
      {\|w_j\|_{L^2(\Omega)}^2}.
    \label{eq:projected_p}
\end{equation}
This provides an adjoint-free, closed-form directional indicator predicting the first-order change in $\mathcal{Q}$ along the candidate direction $\phi_j$. It is, however, not intended to provide an accurate approximation of $\mathcal{Q}(V)-\mathcal{Q}(\hat V^{(k)})$, but rather as an easy-to-compute score for selecting the next active basis function.

\subsection{Residual-decay indicator}
\label{sec:r_decay}
We now focus on the choice $\mathcal{Q}(V):=\frac{1}{2}\|\mathcal{L}(V)\|^2_{L^2(\Omega)}$, which leads to the following indicator.
\begin{definition}[Residual-decay indicator]
\label{def:awr}
Let $\hat V^{(k)}\in\mathcal V_{n_k}$, $\mathcal L_k:=D\mathcal L(\hat V^{(k)})$ and residual $R(\hat V^{(k)}):=\mathcal L(\hat V^{(k)})$.
For each candidate $\phi_j, \; j\in\Lambda^{(k)}$, the residual-decay indicator is given by
\begin{equation}
e_{j,\mathrm{PDE}}^{(k)}:=\frac{\big|\langle R(\hat V^{(k)}),\,\mathcal L_k\phi_j\rangle_{L^2(\Omega)}\big|}{\|\mathcal L_k\phi_j\|_{L^2(\Omega)}}.
\label{eq:basis_indicator}
\end{equation}
 We set $e^{(k)}_{j,\mathrm{PDE}} = 0$ whenever
 $\|\mathcal{L}_k\phi_j\|_{L^2(\Omega)} = 0$.
\end{definition}
To interpret the indicator, consider for each $j\in\Lambda^{(k)}$ the single-basis enrichment $\hat V^{(k)}+\alpha\phi_j$, $\alpha\in\mathbb R$, and the quadratic model of the residual functional obtained by linearizing $\mathcal L$ at $\hat V^{(k)}$,
\begin{equation}
m_j(\alpha):=\frac{1}{2}\bigl\|R\bigl(\hat{V}^{(k)}\bigr)+\alpha\, \mathcal{L}_k \phi_j\bigr\|_{L^2(\Omega)}^2,\qquad\alpha\in\mathbb R.
\label{eq:Jkp1}
\end{equation}
The next proposition shows that $e_{j,\mathrm{PDE}}^{(k)}$ is, up to sign and a square, the projected DWR indicator \eqref{eq:projected_p}, and that it measures the largest reduction of the model $m_j$ achievable along $\phi_j$.

\begin{proposition}[Residual-decay indicator as predicted reduction]
\label{prop:rdecay}
Let $\mathcal{Q}(V)=\frac{1}{2}\|\mathcal{L}(V)\|^2_{L^2(\Omega)}$ and let $j\in\Lambda^{(k)}$ with $\mathcal L_k\phi_j\neq0$. Then the following hold.
\begin{enumerate}[label=(\roman*),leftmargin=*]
\item $\mathcal E_j^{\mathrm{proj}}=-\bigl(e_{j,\mathrm{PDE}}^{(k)}\bigr)^2$.
\item $\mathcal Q(\hat V^{(k)}+\alpha\phi_j)=m_j(\alpha)+o(|\alpha|)$ as $\alpha\to0$.
\item The model $m_j$ has the unique minimizer
\begin{equation*}
\alpha_j^*=-\frac{\bigl\langle R\bigl(\hat{V}^{(k)}\bigr), \mathcal{L}_k \phi_j\bigr\rangle_{L^2(\Omega)}}{\left\|\mathcal{L}_k \phi_j\right\|_{L^2(\Omega)}^2},
\end{equation*}
and $\mathcal Q\bigl(\hat V^{(k)}\bigr)-m_j(\alpha_j^*)=\frac12\bigl(e_{j,\mathrm{PDE}}^{(k)}\bigr)^2$.
\item $e_{j,\mathrm{PDE}}^{(k)}=\|\Pi_jR(\hat V^{(k)})\|_{L^2(\Omega)}$, where $\Pi_j$ denotes the $L^2(\Omega)$-orthogonal projection onto $\operatorname{span}\{\mathcal L_k\phi_j\}$. In particular, $e_{j,\mathrm{PDE}}^{(k)}$ is unchanged if $\phi_j$ is replaced by $s\phi_j$ with $s\neq0$.
\end{enumerate}
\end{proposition}
\begin{proof}
Write $R:=R(\hat V^{(k)})$ and $w_j:=\mathcal L_k\phi_j$, and let $\langle\cdot,\cdot\rangle$ and $\|\cdot\|$ denote the inner product and norm of $L^2(\Omega)$.

(i) By the chain rule, $\langle\mathcal Q'(\hat V^{(k)}),\psi\rangle_{\mathcal V^*,\mathcal V}=\langle R,\mathcal L_k\psi\rangle$ for all $\psi\in\mathcal V$. Taking $\psi=\phi_j$ in \eqref{eq:projected_p} gives $\mathcal E_j^{\mathrm{proj}}=-\langle R,w_j\rangle^2/\|w_j\|^2$.

(ii) With the remainder $r_{\mathcal L}$ introduced before Lemma~\ref{lem:dwr}, $\mathcal L(\hat V^{(k)}+\alpha\phi_j)=R+\alpha w_j+r_{\mathcal L}(\alpha\phi_j)$, hence
\begin{equation*}
\mathcal Q\bigl(\hat V^{(k)}+\alpha\phi_j\bigr)-m_j(\alpha)=\bigl\langle R+\alpha w_j,\,r_{\mathcal L}(\alpha\phi_j)\bigr\rangle+\tfrac12\|r_{\mathcal L}(\alpha\phi_j)\|^2 .
\end{equation*}
Since $\|r_{\mathcal L}(\alpha\phi_j)\|=o(|\alpha|\,\|\phi_j\|_{\mathcal V})$, the right-hand side is $o(|\alpha|)$.

(iii) Expanding the square, $m_j(\alpha)=\frac12\|R\|^2+\alpha\langle R,w_j\rangle+\frac12\alpha^2\|w_j\|^2$, a strictly convex quadratic since $w_j\neq0$. Its derivative vanishes only at $\alpha_j^*$, and substitution gives $m_j(\alpha_j^*)=\frac12\|R\|^2-\frac12\langle R,w_j\rangle^2/\|w_j\|^2$. The claim follows from $m_j(0)=\mathcal Q(\hat V^{(k)})$.

(iv) We have $\Pi_jR=\langle R,w_j\rangle\|w_j\|^{-2}w_j$, whose norm is $e_{j,\mathrm{PDE}}^{(k)}$. Since $\mathcal L_k$ is linear, replacing $\phi_j$ by $s\phi_j$ replaces $w_j$ by $sw_j$, which leaves $\operatorname{span}\{w_j\}$ and hence $\Pi_j$ unchanged. \qed
\end{proof}
By Proposition~\ref{prop:rdecay}(iii), selecting $\phi_j$ with the largest value of $e_{j,\mathrm{PDE}}^{(k)}$ maximizes the reduction of the residual functional $\mathcal{Q}$ predicted by the linearized model $m_j$. Part~(ii) only shows that $m_j$ is accurate to first order in $\alpha$; the indicator is therefore a ranking criterion rather than a guaranteed decrease of $\mathcal Q$.

\paragraph{Greedy selection of basis}
At the $k$-th iteration, among all candidates  $\{\phi_j: j\in\Lambda^{(k)}\}$,
the next basis function is then chosen as
$$
j^*:=\operatorname{argmax}_{j \in\Lambda^{(k)}}\; e_{j,\mathrm{PDE}}^{(k)}.
$$
The active and candidate index sets are then updated as
\begin{equation*}
\Gamma^{(k+1)}=\Gamma^{(k)}\cup \{j^*\},\quad\Lambda^{(k+1)}=\Lambda^{(k)}\setminus\{j^*\}.
\end{equation*}
In practice, after adding a new basis function $\phi_{j^*}$, we re-optimize all coefficients over $\mathcal{V}_{n_k+1}$, so that $\hat{V}^{(k+1)}$ is not exactly $\hat{V}^{(k)}+\alpha_{j^*}^* \phi_{j^*}$. Nevertheless, by Proposition~\ref{prop:rdecay}(iv), $e_{j,\mathrm{PDE}}^{(k)}$ is the $L^2(\Omega)$-norm of the projection of the current residual onto the one-dimensional space spanned by $\mathcal L_k\phi_j$. It therefore selects the candidate whose linearized image $\mathcal L_k\phi_j$ captures the largest component of the current residual, providing an effectively computable indicator for basis selection.

\subsection{The stationary HJB equation}
\label{sec:HJB_indicator}
We apply the residual-decay indicator settings to the stationary HJB equation \eqref{eq:stationary} associated with the infinite-horizon problem \eqref{eq:infinite_horizon}. Since the basis functions depend only on the state, and the finite-horizon value function at $t=0$ approximates the stationary one for large $T$ (Section~\ref{sec:ocp}), we score candidate basis functions with the stationary HJB residual. By the Dynamic Programming Principle, the optimal feedback control is given by
\begin{equation*}
\mathbf{u}^*(\mathbf{x})=-\frac{1}{2 \beta} \mathbf{q}(\mathbf{x})^\top \nabla V(\mathbf{x}).
\end{equation*}
Substituting $\mathbf{u}^*$ into  \eqref{eq:stationary} gives the stationary HJB operator $\mathcal{L}_{\mathrm{HJB}}: \mathcal{V}\rightarrow L^{2}(\Omega)$:
\begin{equation}
(\mathcal{L}_{\mathrm{HJB}}V)(\mathbf{x}):=-\frac{1}{4 \beta} \nabla V(\mathbf{x})^\top \mathbf{q}(\mathbf{x}) \mathbf{q}(\mathbf{x})^\top \nabla V(\mathbf{x})+\nabla V(\mathbf{x})^\top \mathbf{f}(\mathbf{x})+\ell(\mathbf{x}),
\label{eq:HJBoperator}
\end{equation}
for any $\mathbf{x} \in \Omega$.
Let $\hat V^{(k)}\in\mathcal V_{n_k}$ denote the $k$-th approximation of the value function, with index set $\Gamma^{(k)}$ satisfying $|\Gamma^{(k)}|=n_k$. The associated HJB residual is defined as $R(\hat V^{(k)})(\mathbf{x}):= \big(\mathcal{L}_{\mathrm{HJB}}\hat V^{(k)}\big)(\mathbf{x})$.
Since $\mathcal L_{\mathrm{HJB}}$ is quadratic in $\nabla V$, its expansion around $\hat V^{(k)}$ terminates at second order. This gives the linearized operator explicitly.

\begin{lemma}[Linearization of the HJB operator]
\label{lem:HJBlin}
Let $\hat V^{(k)}\in\mathcal V_{n_k}$ and $\mathbf{a}_k(\mathbf{x}):=\mathbf{f}(\mathbf{x}) -\frac{1}{2\beta}\,\mathbf{q}(\mathbf{x})\mathbf{q}(\mathbf{x})^\top \nabla \hat V^{(k)}(\mathbf{x})$. For every $\psi\in\mathcal V$, $\varepsilon\in\mathbb R$ and $\mathbf x\in\Omega$,
\begin{equation}
\mathcal L_{\mathrm{HJB}}\bigl(\hat V^{(k)}+\varepsilon\psi\bigr)(\mathbf x)
=R\bigl(\hat V^{(k)}\bigr)(\mathbf x)+\varepsilon\,\mathbf a_k(\mathbf x)^\top\nabla\psi(\mathbf x)-\frac{\varepsilon^2}{4\beta}\bigl\|\mathbf q(\mathbf x)^\top\nabla\psi(\mathbf x)\bigr\|_2^2 .
\label{eq:HJBexpansion}
\end{equation}
Consequently, the linearized HJB operator $\mathcal L_{\mathrm{HJB},k}:\mathcal V\to L^2(\Omega)$ at $\hat V^{(k)}$ is
\begin{equation}
(\mathcal L_{\mathrm{HJB},k}\psi)(\mathbf{x}) := D\mathcal{L}_{\mathrm{HJB}}\bigl(\hat V^{(k)}\bigr)[\psi](\mathbf{x}) = \mathbf{a}_k(\mathbf{x})^\top\nabla\psi(\mathbf{x}).
\label{eq:LkHJB}
\end{equation}
\end{lemma}
\begin{proof}
Fix $\mathbf x\in\Omega$ and write $\mathbf p:=\nabla\hat V^{(k)}(\mathbf x)$, $\mathbf g:=\nabla\psi(\mathbf x)$ and $\mathbf B:=\mathbf q(\mathbf x)\mathbf q(\mathbf x)^\top$, which is symmetric. Then
\begin{equation*}
(\mathbf p+\varepsilon\mathbf g)^\top\mathbf B(\mathbf p+\varepsilon\mathbf g)=\mathbf p^\top\mathbf B\mathbf p+2\varepsilon\,\mathbf g^\top\mathbf B\mathbf p+\varepsilon^2\|\mathbf q(\mathbf x)^\top\mathbf g\|_2^2 .
\end{equation*}
Inserting this into \eqref{eq:HJBoperator} and collecting the terms linear in $\varepsilon$ gives $\varepsilon\,\mathbf g^\top\bigl(\mathbf f(\mathbf x)-\frac{1}{2\beta}\mathbf B\mathbf p\bigr)=\varepsilon\,\mathbf a_k(\mathbf x)^\top\mathbf g$, which proves \eqref{eq:HJBexpansion}. Dividing by $\varepsilon$ and letting $\varepsilon\to0$ in \eqref{eq:HJBexpansion} identifies the derivative of $\mathcal L_{\mathrm{HJB}}$ at $\hat V^{(k)}$ in the direction $\psi$, yielding \eqref{eq:LkHJB}. \qed
\end{proof}
By \eqref{eq:HJBexpansion}, the linearized residual $R(\hat V^{(k)})+\alpha\,\mathcal L_{\mathrm{HJB},k}\phi_j$ used in the model \eqref{eq:Jkp1} differs from the exact residual $\mathcal L_{\mathrm{HJB}}(\hat V^{(k)}+\alpha\phi_j)$ only by the term $-\frac{\alpha^2}{4\beta}\|\mathbf q^\top\nabla\phi_j\|_2^2$, which is small when the control acts weakly in the direction of $\nabla\phi_j$.
We now define a stationary HJB residual-decay indicator.
\begin{definition}[Stationary HJB residual-decay indicator]
\label{def:HJB}
For the infinite-horizon problem \eqref{eq:infinite_horizon}, let $\mathcal{L}_{\mathrm{HJB}}$ be given by \eqref{eq:HJBoperator}  and its associated residual $R(\hat V^{(k)})$. For $\hat V^{(k)}\in\mathcal V_{n_k}$, let $\mathbf{a}_k$ be as in Lemma~\ref{lem:HJBlin}. For each candidate basis function $\phi_j$ with index $j\in\Lambda^{(k)}:=\Theta\setminus\Gamma^{(k)}$, the stationary HJB residual-decay indicator is defined by
\begin{equation}
\begin{aligned}
e_{j,\mathrm{HJB}}^{(k)}:=\frac{\left|\int_{\Omega} R\left(\hat{V}^{(k)}\right)(\mathbf{x}) \left(\mathbf{a}_k(\mathbf{x})^\top\nabla \phi_j(\mathbf{x})\right) d \mathbf{x}\right|}{\left(\int_{\Omega}\left|\mathbf{a}_k(\mathbf{x})^\top \nabla \phi_j(\mathbf{x})\right|^2 d \mathbf{x}\right)^{1 / 2}}.
\end{aligned}
\label{eq:indicatorHJB}
\end{equation}
We set $e_{j,\mathrm{HJB}}^{(k)}=0$ whenever the denominator in \eqref{eq:indicatorHJB} vanishes.
\end{definition}
\begin{remark}[Constant basis functions]
\label{rem:constant}
The indicator $e_{j,\mathrm{HJB}}^{(k)}$ depends only on $\nabla \phi_j$ and is therefore ``blind'' to constant basis functions. This is consistent with the fact that the undiscounted stationary HJB equation \eqref{eq:stationary} determines $V$ only up to additive constants: the indicator effectively learns only the shape of $V$, i.e., its gradient. Since the operator is used solely to select basis functions, and not to solve the stationary problem, we handle this by requiring the initial active set $\Gamma^{(1)}$ to contain a constant basis function.
\end{remark}

\subsection{Model-based greedy algorithm}
After introducing the indicator for adaptive basis selection, we now define a validation error used to analyze the approximation error throughout the iterative procedure. Let the validation dataset be $\mathcal{D}_{\mathrm{val}}=\left\{\left(\mathbf{x}_{\mathrm{val},i}, v_{\mathrm{val},i}, \mathbf{g}_{\mathrm{val},i}\right)\right\}_{i=1}^{N_{\mathrm{val}}}$ with $N_{\mathrm{val}}$ points. The relative $L^2$- and $H^1$-validation errors of $\hat{V}^{(k)}$ at the $k$-th iteration are computed as
\begin{equation}
	\label{eq:kk1}
	\begin{split}
		\operatorname{Error}_{L^2}(\hat{V}^{(k)})  &\displaystyle =\left(\frac{\displaystyle\sum\limits_{i=1}^{N_{\mathrm{val}}}\left| \hat{V}^{(k)}(\mathbf{x}_{\mathrm{val},i})-v_{\mathrm{val},i}\right|^2}{\displaystyle\sum\limits_{i=1}^{N_{\mathrm{val}}}\left| v_{\mathrm{val},i}\right|^2}\right)^{\frac12} , \\
		\operatorname{Error}_{H^1}(\hat{V}^{(k)}) &=\left(\frac{\displaystyle\sum\limits_{i=1}^{N_{\mathrm{val}}}\displaystyle\left\|  \nabla\hat{V}^{(k)}(\mathbf{x}_{\mathrm{val},i})  -\mathbf{g}_{\mathrm{val},i}\right\|_2^2 }{\displaystyle\sum\limits_{i=1}^{N_{\mathrm{val}}}\left\|\mathbf{g}_{\mathrm{val},i}\right\|_2^2 }\right)^{\frac12}.
	\end{split}
\end{equation}
With both the selection criterion and error measures established, we are now ready to present the model-based greedy algorithm.
\begin{algorithm}[ht]
\caption{Model-based greedy algorithm}
\begin{algorithmic}[1]
\State \textbf{Input:} Candidate basis set $\Phi_{\Theta}$, validation dataset $\mathcal{D}_{\mathrm{val}}$, iteration budget $k_{\max}$, early stop tolerance $\varepsilon$
\State \textbf{Initialize} active basis index set $\Gamma^{(1)}$ with cardinality $1\leq |\Gamma^{(1)}|\ll P$, and candidate set $\Lambda^{(1)}=\Theta\setminus\Gamma^{(1)}$
\State Solve for coefficients $\mathbf{c}^{(1)}$ by Galerkin methods on $\Phi_{\Gamma^{(1)}}:=\{\phi_j: j\in\Gamma^{(1)}\}$
\State Set $ \hat V^{(1)}(\mathbf{x}) = \Phi_{\Gamma^{(1)}}(\mathbf{x})^\top\mathbf{c}^{(1)}$
 \For{$k=1,2,\dots,k_{\max}$}
  \State Compute $e^{(k)}_{j,\mathrm{PDE}}$ for $j\in\Lambda^{(k)}$
  \If{$\max_{j\in\Lambda^{(k)}}e^{(k)}_{j,\mathrm{PDE}}<\varepsilon$}
    \State \textbf{break}
  \EndIf
  \State $j^*\gets\operatorname{argmax}_{j\in\Lambda^{(k)}}e^{(k)}_{j,\mathrm{PDE}}$
  \State $\Gamma^{(k+1)}=\Gamma^{(k)}\cup\{j^*\}$ and $\Lambda^{(k+1)}=\Lambda^{(k)}\setminus\{j^*\}$
  \State Solve for $\mathbf{c}^{(k+1)}$ on $\Phi_{\Gamma^{(k+1)}}$ and set
    $\hat V^{(k+1)}=\Phi_{\Gamma^{(k+1)}}^\top\mathbf{c}^{(k+1)}$
  \State Compute $\mathrm{Error}_{L^2}(\hat V^{(k+1)})$ and $\mathrm{Error}_{H^1}(\hat V^{(k+1)})$
    on $\mathcal{D}_{\mathrm{val}}$
\EndFor
\end{algorithmic}
\label{algorithm_1}
\end{algorithm}

Algorithm~\ref{algorithm_1} describes the model-based greedy algorithm for solving the PDE \eqref{eq:pde} by Galerkin methods, where the basis selection is guided by the residual-decay indicator. Starting from an initial basis index set, the algorithm iteratively selects the candidate with the largest value of the indicator. At each iteration, the Galerkin system is solved in the updated basis subspace, and relative validation errors are computed. The procedure terminates either when the prescribed iteration budget $k_{\max}$ is reached or when all remaining indicator values fall below the tolerance $\varepsilon$, i.e., $\max _{j \in \Lambda^{(k)}}e_{j,\mathrm{PDE}}^{(k)}<\varepsilon$.

\begin{remark}[Computational cost]\label{rem:cost}
At each iteration of Algorithm~\ref{algorithm_1} we compute
$e^{(k)}_{j,\mathrm{PDE}}$ for all $j\in\Lambda^{(k)}$ and then re-solve the
Galerkin system on the enlarged subspace. Using quadrature with $N_q$ nodes,
the computation of all
candidate indicators costs $O\bigl(N_q\,|\Lambda^{(k)}|\,d\bigr)$ operations.
 This step uses
the current residual values and projects them against the linearized candidate
directions. The Galerkin problem is then re-solved only on the active subspace
of dimension $n_k\ll P$.

\noindent The cumulative cost of the linear-algebra component over $n$ greedy iterations is
$\sum_{k=1}^{n}O(n_k^3)$, in contrast, solving the full non-adaptive Galerkin model with all $P$ basis
functions has cost typically of order
$O(P^3)$. For large $P$, the reduced Galerkin solves can provide substantial
computational savings. 

\noindent Finally, as discussed in Section~\ref{sec:data}, when training data are available, $e_{j,\mathrm{PDE}}^{(k)}$ can be approximated directly from data, further reducing the cost of the selection step and making the method particularly attractive in data-driven settings.

\end{remark}

\section{Data-driven adaptive basis and data selection}
\label{sec:data}
Given a training dataset $\mathcal{D}=\left\{\left(\mathbf{x}_i, v_i, \mathbf{g}_i\right)\right\}_{i=1}^N$, we introduce an analogous partition of the data in \begin{equation*}
\mathcal{I}^{(k)} \cup \mathcal{J}^{(k)}=\{1, \cdots, N\}, \quad \mathcal{I}^{(k)} \cap \mathcal{J}^{(k)}=\emptyset,
\end{equation*}
where $\mathcal{I}^{(k)} \subset\{1, \cdots, N\}$ indexes the active training data, \[\mathcal{D}_{\mathcal{I}^{(k)}}=\left\{\left(\mathbf{x}_i, v_i, \mathbf{g}_i\right)\right\}_{i \in \mathcal{I}^{(k)}},\] used to fit $\hat{V}^{(k)}$, with cardinality $m_k:=\left|\mathcal{I}^{(k)}\right|$. The set $\mathcal{J}^{(k)}=\{1, \cdots, N\} \backslash \mathcal{I}^{(k)}$ indexes the remaining data $\mathcal{D}_{\mathcal{J}^{(k)}}=\mathcal{D} \setminus \mathcal{D}_{\mathcal{I}^{(k)}}$, which are used for residual estimation and adaptive data selection.

\noindent We introduce data-driven criteria for adaptive basis and data selection within the supervised learning framework presented in Section~\ref{sec:sp}. At each iteration, the most informative data points and basis functions are selected according to their predicted contribution to reducing the supervised loss.

\subsection{Supervised-loss basis indicator}
To quantify the potential improvement obtained by adding a single candidate basis function, we use the predicted one-step decrease of the gradient-augmented
supervised loss. The construction is analogous to the single-direction projected indicator in Section~\ref{sec:single}. In the model-based setting, the
PDE residual is tested against the linearized image $\mathcal L_k\phi_j$ of a candidate basis direction. In the supervised setting, the PDE residual is replaced
by the empirical residual on the active training set, and the linearized image is
replaced by the corresponding feature vector generated by $\phi_j$ and
$\nabla\phi_j$. Thus, in both cases, each candidate basis function is assessed by
the improvement it can produce along its own direction.

\noindent For the active dataset $\mathcal D_{\mathcal I^{(k)}}$, let $\mathbf h^{(k)}$ denote the weighted residual vector associated with the
current approximation on the active training set $\mathcal I^{(k)}$, collecting both
function value and gradient residuals:
\begin{equation*}
\mathbf h^{(k)}=\left[\begin{array}{c}
\{\sqrt{\omega_1\mu_i}\,(v_i-\hat V^{(k)}(\mathbf x_i))\}_{i\in\mathcal I^{(k)}}\\
\{\sqrt{\omega_2\mu_i}\,(\mathbf g_i-\nabla\hat V^{(k)}(\mathbf x_i))\}_{i\in\mathcal I^{(k)}}
\end{array}\right]\in\mathbb R^{(1+d)m_k},
\end{equation*}
where $\mu_i>0$ denotes the sampling weight (for uniformly sampled data one may take $\mu_i=1/m_k$ or $\mu_i=1$), and $\omega_1,\omega_2>0$ weight the function-value and gradient residuals, respectively. Similarly, for each candidate
$\phi_j$, $j\in\Lambda^{(k)}$, define the weighted feature vector
\begin{equation*}
\mathbf z_j^{(k)}=\left[\begin{array}{c}
\{\sqrt{\omega_1\mu_i}\,\phi_j(\mathbf x_i)\}_{i\in\mathcal I^{(k)}}\\
\{\sqrt{\omega_2\mu_i}\,\nabla\phi_j(\mathbf x_i)\}_{i\in\mathcal I^{(k)}}
\end{array}\right]\in\mathbb R^{(1+d)m_k}.
\end{equation*}
For each candidate basis function $\phi_j$, $j\in\Lambda^{(k)}$, we consider the one-direction enrichment $\hat V^{(k)}+\alpha\phi_j$, $\alpha\in\mathbb R$. The supervised loss on the active training set, restricted to this enrichment, is
\begin{equation*}
\mathcal{Q}_{\mathrm{SL},j}^{(k)}(\alpha):=\frac12\bigl\|\mathbf h^{(k)}-\alpha\,\mathbf z_j^{(k)}\bigr\|_2^2,\qquad\alpha\in\mathbb R,
\end{equation*}
so that $\mathcal{Q}_{\mathrm{SL},j}^{(k)}(0)$ is the current supervised loss. We define the supervised-loss basis indicator
\begin{align}
 e_{j,\mathrm{SL}}^{(k)}
&:=
\frac{
\left|
\left(\mathbf z_j^{(k)}\right)^\top\mathbf h^{(k)}
\right|
}{
\left\|\mathbf z_j^{(k)}\right\|_2
} \nonumber\\
&=
 \left(
 \sum_{i\in\mathcal I^{(k)}}\mu_i
 \left[
 \omega_1|\phi_j(\mathbf x_i)|^2
 +
 \omega_2\|\nabla\phi_j(\mathbf x_i)\|_2^2
 \right]
 \right)^{-1/2}\notag\\
&\quad\times
 \Bigg|
 \sum_{i\in\mathcal I^{(k)}}\mu_i
 \Big[
 \omega_1
 \left(v_i-\hat V^{(k)}(\mathbf x_i)\right)\phi_j(\mathbf x_i)\notag\\
&\qquad\qquad\qquad +
 \omega_2
 \left(\mathbf g_i-\nabla\hat V^{(k)}(\mathbf x_i)\right)^\top
 \nabla\phi_j(\mathbf x_i)
 \Big]
 \Bigg|,
\label{eq:supervised_loss_indicator}
\end{align}
with $e_{j,\mathrm{SL}}^{(k)}:=0$ whenever $\mathbf z_j^{(k)}=\mathbf 0$. It is the data-driven counterpart of the residual-decay indicator, as the following proposition shows.

\noindent To compare with OMP, let $\mathbf r^{(k)}$ and $\mathbf a_j$ denote the vectors $\mathbf h^{(k)}$ and $\mathbf z_j^{(k)}$ computed with $\omega_1=\omega_2=\mu_i=1$, that is, the unweighted residual vector and the unweighted gradient-augmented column of $\phi_j$ on the active data. Applied to the dictionary $\{\mathbf a_j\}_{j\in\Lambda^{(k)}}$ with residual $\mathbf r^{(k)}$, the selection step of OMP chooses
\begin{equation}
j^*_{\mathrm{OMP}}\in\operatorname{argmax}_{j\in\Lambda^{(k)}}\frac{|\mathbf a_j^\top\mathbf r^{(k)}|}{\|\mathbf a_j\|_2}.
\label{eq:omp_rule}
\end{equation}

\begin{proposition}[Supervised-loss indicator]
\label{prop:sl}
Let $j\in\Lambda^{(k)}$ with $\mathbf z_j^{(k)}\neq\mathbf 0$. Then the following hold.
\begin{enumerate}[label=(\roman*),leftmargin=*]
\item $\mathcal{Q}_{\mathrm{SL},j}^{(k)}$ has the unique minimizer $\alpha_j^*=(\mathbf z_j^{(k)})^\top\mathbf h^{(k)}/\|\mathbf z_j^{(k)}\|_2^2$, and $\mathcal{Q}_{\mathrm{SL},j}^{(k)}(0)-\mathcal{Q}_{\mathrm{SL},j}^{(k)}(\alpha_j^*)=\frac12\bigl(e_{j,\mathrm{SL}}^{(k)}\bigr)^2$.
\item $e_{j,\mathrm{SL}}^{(k)}$ is the Euclidean norm of the orthogonal projection of $\mathbf h^{(k)}$ onto $\operatorname{span}\{\mathbf z_j^{(k)}\}$, and it is unchanged if $\phi_j$ is replaced by $s\phi_j$ with $s\neq0$.
\item If $\omega_1=\omega_2$ and the weights $\mu_i$, $i\in\mathcal I^{(k)}$, are equal, then $\operatorname{argmax}_{j\in\Lambda^{(k)}}e_{j,\mathrm{SL}}^{(k)}$ coincides with the OMP selection rule \eqref{eq:omp_rule}.
\end{enumerate}
\end{proposition}
\begin{proof}
Parts (i) and (ii) follow as in Proposition~\ref{prop:rdecay}(iii)--(iv), with $L^2(\Omega)$ replaced by $\mathbb R^{(1+d)m_k}$ with the Euclidean inner product, and $R(\hat V^{(k)})$ and $\mathcal L_k\phi_j$ replaced by $\mathbf h^{(k)}$ and $-\mathbf z_j^{(k)}$; scale invariance holds because $\mathbf z_j^{(k)}$ is linear in $\phi_j$. For (iii), if $\omega_1=\omega_2=\omega$ and $\mu_i=\mu$, then $\mathbf h^{(k)}=\sqrt{\omega\mu}\,\mathbf r^{(k)}$ and $\mathbf z_j^{(k)}=\sqrt{\omega\mu}\,\mathbf a_j$. Hence $e_{j,\mathrm{SL}}^{(k)}=\sqrt{\omega\mu}\,|\mathbf a_j^\top\mathbf r^{(k)}|/\|\mathbf a_j\|_2$, and the common factor $\sqrt{\omega\mu}>0$ does not affect the maximizer in \eqref{eq:omp_rule}. \qed
\end{proof}
By Proposition~\ref{prop:sl}(i), selecting the candidate with the largest $e_{j,\mathrm{SL}}^{(k)}$ selects the candidate that yields the largest one-step decrease of the supervised loss on the active data. Unlike Proposition~\ref{prop:rdecay}(iii), this decrease is exact, since the supervised loss is quadratic in the coefficients. By Proposition~\ref{prop:sl}(iii), the basis step of the data-driven algorithm differs from OMP only through the weights $\omega_1,\omega_2$ and the LASSO refit; the additional ingredient is the greedy selection of the training data in Section~\ref{sec:data_score}. The next basis function is chosen as
\begin{equation*}
j^*=\operatorname{argmax}_{j\in\Lambda^{(k)}} e_{j,\mathrm{SL}}^{(k)},
\end{equation*}
and the index sets are updated by
\begin{equation*}
\Gamma ^{(k+1)}=\Gamma^{(k)} \cup \{j^*\}, \quad \Lambda^{(k+1)}=\Lambda^{(k)} \backslash \{j^*\}.
\end{equation*}
After each enrichment step, all coefficients are re-estimated on the enlarged active basis set. Alternatively, one may assemble the candidate columns $\mathbf Z^{(k)}\in\mathbb{R}^{(1+d)  m_k \times |\Lambda^{(k)}|}$, and solve a regularized least-squares problem over the  whole candidate space, in that case, the resulting coefficients should be interpreted as a proxy for multi-direction correction rather than as the exact one-step loss reduction with any individual candidate basis function.

\begin{remark}[Data-driven residual-decay indicator]
\label{rem:empiricalPDE}
The residual-decay indicator of Definition~\ref{def:awr} also has a data-driven version. Replacing the $L^2(\Omega)$ inner product in \eqref{eq:basis_indicator} by the $\mu$-weighted discrete inner product on the active points $\{\mathbf{x}_i\}_{i \in \mathcal{I}^{(k)}}$ gives
\begin{equation}
\tilde{e}_{j,\mathrm{PDE}}^{(k)}:=\frac{\left|\sum_{i \in \mathcal{I}^{(k)}} \mu_i R(\hat{V}^{(k)})(\mathbf{x}_i)\,\mathcal{L}_k \phi_j(\mathbf{x}_i)\right|}{\left(\sum_{i \in \mathcal{I}^{(k)}} \mu_i\left|\mathcal{L}_k \phi_j(\mathbf{x}_i)\right|^2\right)^{1 / 2}}.
\label{eq:empiricalPDEResidue}
\end{equation}
The analogues of Proposition~\ref{prop:rdecay}(iii)--(iv) for $\tilde{e}_{j,\mathrm{PDE}}^{(k)}$, with this inner product in place of that of $L^2(\Omega)$, hold by the same proof. This provides a practical way to use the DWR-type indicators in high-dimensional, data-driven settings.
\end{remark}

\subsection{Data score}
\label{sec:data_score}
In contrast to the basis-selection indicator, which estimates the loss decrease
obtained by adding a candidate direction, the data score measures the current
local contribution of each unselected sample to the supervised residual.
For each candidate training data point $(\mathbf{x}_i,v_i,\mathbf{g}_i) \in \mathcal{D}_{\mathcal{J}^{(k)}}$, we define the gradient-augmented data score in the same spirit as the supervised-loss basis indicator:
\begin{equation}
\operatorname{score}_i^{(k)}=\omega_1\left|v_i-\hat{V}^{(k)}\left(\mathbf{x}_i\right)\right|^2+\omega_2\left\|\mathbf g_i-\nabla \hat{V}^{(k)}\left(\mathbf{x}_i\right)\right\|_2^2.
\label{eq:data_score}
\end{equation}
We then select the next data point as
\begin{equation*}
i^*=\operatorname{argmax}_{i \in \mathcal{J}^{(k)}} \operatorname{score}_i^{(k)},
\end{equation*}
and update the data index sets
\begin{equation*}
\mathcal{I} ^{(k+1)}=\mathcal{I}^{(k)} \cup \{i^*\}, \quad \mathcal{J}^{(k+1)}=\mathcal{J}^{(k)} \backslash\{i^*\}.
\end{equation*}

\subsection{Data-driven greedy algorithm}

\begin{algorithm}[ht]
\caption{Data-driven greedy algorithm}
\begin{algorithmic}[1]
\State \textbf{Input:} Candidate basis set $\Phi_{\Theta}$, training dataset $\mathcal{D}$, validation dataset $\mathcal{D}_{\mathrm{val}}$, iteration budget $k_{\max}$, thresholds $\varepsilon_{\mathrm{basis}}$ and $\varepsilon_{\mathrm{data}}$
\State \textbf{Initialize} active basis set $\Gamma^{(1)}$, candidate basis set $\Lambda^{(1)}=\Theta\setminus\Gamma^{(1)}$, active dataset $\mathcal{I}^{(1)}$, inactive dataset $\mathcal J^{(1)}=\{1,\dots,N\}\setminus \mathcal{I}^{(1)}$, and \texttt{basis\_stop} = \texttt{FALSE}, \texttt{data\_stop} = \texttt{FALSE}
\State Solve for $\mathbf{c}^{(1)}$ by supervised learning on $\mathcal{D}_{\mathcal{I}^{(1)}}$ and $\Phi_{\Gamma^{(1)}}$, set $\hat V^{(1)}(\mathbf{x}) = \Phi_{\Gamma^{(1)}}(\mathbf{x})^\top\mathbf{c}^{(1)}$
\For{$k=1,2,\ldots,k_{\max}$}
\State Set $\mathcal{I}^{(k+1)}\gets\mathcal{I}^{(k)}$, $\mathcal{J}^{(k+1)}\gets\mathcal{J}^{(k)}$, $\Gamma^{(k+1)}\gets\Gamma^{(k)}$, $\Lambda^{(k+1)}\gets\Lambda^{(k)}$, and $\hat V^{(k+1/2)}\gets\hat V^{(k)}$
\If{\texttt{data\_stop} = \texttt{FALSE} and $\mathcal J^{(k)}\ne\emptyset$}
    \State Compute $\operatorname{score}_i^{(k)}$ for $i \in \mathcal{J}^{(k)}$
    \If{$\max_{i\in\mathcal J^{(k)}} \operatorname{score}_i^{(k)} \geq \varepsilon_{\mathrm{data}}$}
        \State $i^*\gets\operatorname{argmax}_{i \in \mathcal{J}^{(k)}} \operatorname{score}_i^{(k)}$
        \State Update $\mathcal{I}^{(k+1)}=\mathcal{I}^{(k)} \cup \{i^*\}$ and $\mathcal{J}^{(k+1)}=\mathcal{J}^{(k)} \backslash\{i^*\}$
        \State Refit on $\mathcal{D}_{\mathcal{I}^{(k+1)}}$ and $\Phi_{\Gamma^{(k)}}$ to obtain $\mathbf{c}^{(k+1/2)}$ and $\hat V^{(k+1/2)}$
    \Else
        \State \texttt{data\_stop} $\gets$ \texttt{TRUE}
    \EndIf
\EndIf
\State Set $\hat V^{(k+1)}\gets\hat V^{(k+1/2)}$
\If{\texttt{basis\_stop} = \texttt{FALSE} and $\Lambda^{(k)}\ne\emptyset$}
    \State Compute $e_j^{(k)}$ for $j \in \Lambda^{(k)}$ using $\hat V^{(k+1/2)}$
    \If{$\max_{j\in\Lambda^{(k)}} e_j^{(k)} \geq \varepsilon_{\mathrm{basis}}$}
        \State $j^*\gets\operatorname{argmax}_{j \in \Lambda^{(k)}} e_j^{(k)}$
        
        \State Update $\Gamma^{(k+1)} = \Gamma^{(k)} \cup \{j^\ast\}$ and $\Lambda^{(k+1)}=\Lambda^{(k)} \backslash \{j^*\}$
        \State Refit on $\mathcal{D}_{\mathcal{I}^{(k+1)}}$ and $\Phi_{\Gamma^{(k+1)}}$ to obtain $\mathbf{c}^{(k+1)}$ and $\hat V^{(k+1)}$
    \Else
        \State \texttt{basis\_stop} $\gets$ \texttt{TRUE}
    \EndIf
\EndIf
\State Compute $\operatorname{Error}_{L^2}(\hat{V}^{(k+1)})$ and $\operatorname{Error}_{H^1}(\hat{V}^{(k+1)})$ on $\mathcal{D}_{\mathrm{val}}$
\If{\texttt{basis\_stop} = \texttt{TRUE} and \texttt{data\_stop} = \texttt{TRUE}}
    \State \textbf{break}
\EndIf
\EndFor
\end{algorithmic}
\label{algorithm_2}
\end{algorithm}
Algorithm~\ref{algorithm_2} extends the greedy strategy to the supervised learning framework with gradient-augmented data. At each iteration, data scores and basis indicators are evaluated and the relevant coefficients are refitted after any enrichment. Index sets and iterates are carried over unchanged when a step adds nothing, so that the intermediate iterate $\hat V^{(k+1/2)}$ (after the data step) and $\hat V^{(k+1)}$ (after the basis step) are defined at every iteration. The indicator $e_j^{(k)}$ denotes either the supervised-loss basis indicator $e_{j,\mathrm{SL}}^{(k)}$ or the empirical residual-decay indicator $\tilde{e}_{j,\mathrm{PDE}}^{(k)}$. As in Algorithm~\ref{algorithm_1}, separate stopping criteria are imposed for data and basis selection, yielding further computational savings while preserving the desired approximation accuracy.

\section{Numerical experiments}
\label{sec:numerical}
In this section, we present numerical experiments to assess the performance of the greedy selection strategies for basis and training samples introduced above. Section~\ref{sec:testcos} considers benchmark PDEs with closed-form solutions, while Sections~\ref{sec:test1d}--\ref{sec:test6d} consider optimal control problems for which the value functions do not admit closed-form representations. The \texttt{MATLAB} code and data reproducing all numerical experiments in this paper is publicly available at: \url{https://github.com/AmberYuyangHuang/Greedy-algorithm-for-HJB}.

\subsection{A benchmark for basis selection}
\label{sec:testcos}
We first consider the one-dimensional linear PDE
\begin{equation}
\label{eq:1d}
\gamma V(x) = r(x) + b(x) V'(x),\quad x\in\Omega = [-\pi, \pi],
\end{equation}
with $\gamma = 0.1$, $b(x) = 0.05\sin(x)$, and $r$ chosen so that the exact solution is $V(x)=\cos(x)+\cos(10x)$. This is the affine case $\mathcal{L}(V)=\mathcal{A}V-r$ with $(\mathcal{A}\phi)(x)=\gamma \phi(x)-b(x) \phi'(x)$, so the Galerkin system \eqref{eq:galerkin} is linear and $\mathcal{L}_k=\mathcal{A}$ for every $k$. We use the $P=101$ trigonometric basis functions
\begin{equation*}
\Phi_{\Theta}(x) =\frac{1}{\sqrt{\pi}} \left[ \frac{1}{\sqrt{2}}, \cos(x), \sin(x), \dots, \cos(50x), \sin(50x) \right]^\top,
\end{equation*}
with $\Theta=\{1,\dots,101\}$. In the data-driven setting, the training data are $N=60$ equispaced points in $\Omega$, and the validation set consists of $N_{\mathrm{val}}=100$ equispaced points shifted by half a training step.

Since $V$ is supported on only two dictionary elements, $\phi_2=\pi^{-1/2}\cos(x)$ and $\phi_{20}=\pi^{-1/2}\cos(10x)$, this benchmark tests basis identification. Table~\ref{table1} reports the selected indices and the relative $L^2$-validation errors for Algorithm~\ref{algorithm_1} with $e^{(k)}_{j,\mathrm{PDE}}$ and Algorithm~\ref{algorithm_2} with $e^{(k)}_{j,\mathrm{SL}}$, for the initializations $\Gamma^{(1)}=\{2\}$, $\{20\}$ and the constant term $\{1\}$. In all cases, both indicators identify the exact support $\{2,20\}$ within at most two iterations, after which the validation error drops to machine precision.

\begin{table}[htbp]
  \centering
  \footnotesize
  \begin{tabular}{llcccc}
    \toprule
    & & \multicolumn{2}{c}{selected index $j^*$}
      & \multicolumn{2}{c}{$\operatorname{Error}_{L^2}(\hat V^{(k)})$} \\
    \cmidrule(lr){3-4}\cmidrule(lr){5-6}
    $\Gamma^{(1)}$ & indicator & $k=1$ & $k=2$ & $k=1$ & $k=2$ \\
    \midrule
    %$\{2\}$  & $e^{(k)}_{j,\mathrm{PDE}}$ & 20 & 16  & $1.43\times 10^{-15}$ & $1.45\times 10^{-15}$ \\  
    %$\{2\}$  & $e^{(k)}_{j,\mathrm{SL}}$  & 20 & 98  & $1.12\times 10^{-16}$ & $6.59\times 10^{-16}$    \\  
    %$\{20\}$ & $e^{(k)}_{j,\mathrm{PDE}}$ & 2  & 16  &  $1.44\times 10^{-15}$ & $1.46\times 10^{-15} $  \\  
    %$\{20\}$ & $e^{(k)}_{j,\mathrm{SL}}$  & 2  & 98  & $1.03\times 10^{-16}$ & $6.47\times 10^{-16} $ \\ 
    $\{2\}$  & $e^{(k)}_{\cdot,\mathrm{PDE}}$ & 20 & --  & $1.43\times 10^{-15}$ & -- \\ 
    $\{2\}$  & $e^{(k)}_{\cdot,\mathrm{SL}}$  & 20 & --  & $1.12\times 10^{-16}$ & --    \\  
    $\{20\}$ & $e^{(k)}_{\cdot,\mathrm{PDE}}$ & 2  & --  &  $1.44\times 10^{-15}$ &--  \\  
    $\{20\}$ & $e^{(k)}_{\cdot,\mathrm{SL}}$  & 2  & --  & $1.03\times 10^{-16}$ & -- \\  
    
    $\{1\}$  & $e^{(k)}_{\cdot,\mathrm{PDE}}$ & 20 & 2   & $7.38\times 10^{-1}$ & $1.50\times 10^{-15} $  \\  
    $\{1\}$  & $e^{(k)}_{\cdot,\mathrm{SL}}$  & 20 & 2   & $7.02\times 10^{-1}$  & $1.19\times 10^{-16}$  \\  
    \bottomrule
  \end{tabular}
  \caption{Selection paths for the one-dimensional benchmark of
  \eqref{eq:1d}. For each initialization $\Gamma^{(1)}$ we report the
  index $j^*$ selected at iterations $k=1,2$ and the relative
  $L^{2}$-validation error \eqref{eq:kk1} after each refit. The exact
  support of $V$ is $\{2,20\}$. Both the model-based indicator
  $e^{(k)}_{\cdot,\mathrm{PDE}}$ and the supervised-loss indicator
  $e^{(k)}_{\cdot,\mathrm{SL}}$ recover the exact support regardless of
  initialization.}
    \label{table1}
\end{table}

Without prior knowledge of the solution, the initial index can also be chosen by the indicator itself: setting $\hat V^{(0)}=0$, so that $R(\hat V^{(0)})=-r$, and taking $\Gamma^{(1)}=\{j^*\}$ with $j^*\in\operatorname{argmax}_{j\in\Theta}e^{(0)}_{j,\mathrm{PDE}}$ selects $j^*=20$, the dominant component $\cos(10x)$.

\subsection{One-dimensional optimal control problem}
\label{sec:test1d}
\begin{figure}[ht]
\centering
\captionsetup[subfigure]{aboveskip=-1pt,belowskip=-1pt}
  \begin{subfigure}[t]{0.48\textwidth}
\includegraphics[width=\textwidth]{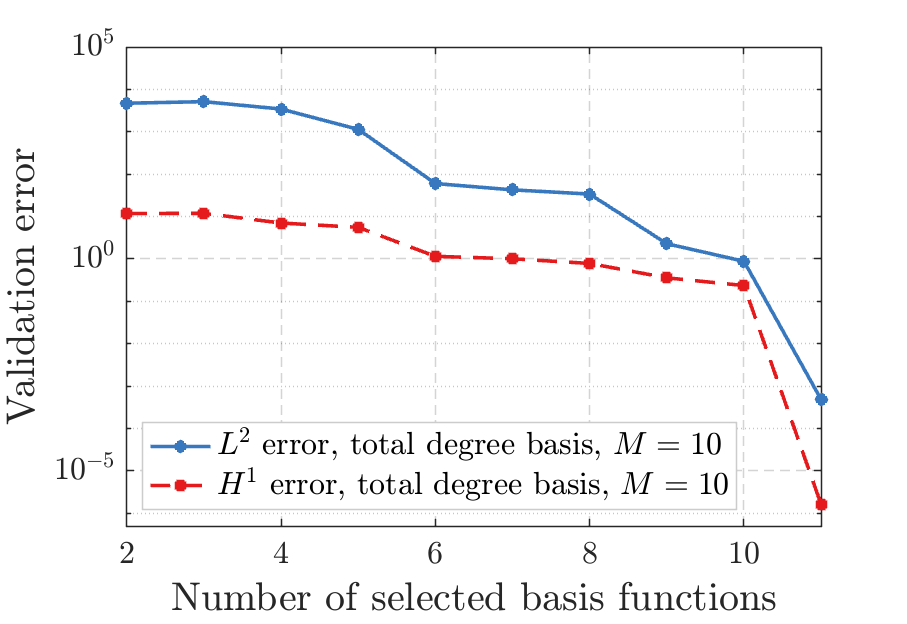}
  \subcaption{Random basis selection}
  \end{subfigure}
\begin{subfigure}[t]
  {0.48\textwidth}
\includegraphics[width=\textwidth]{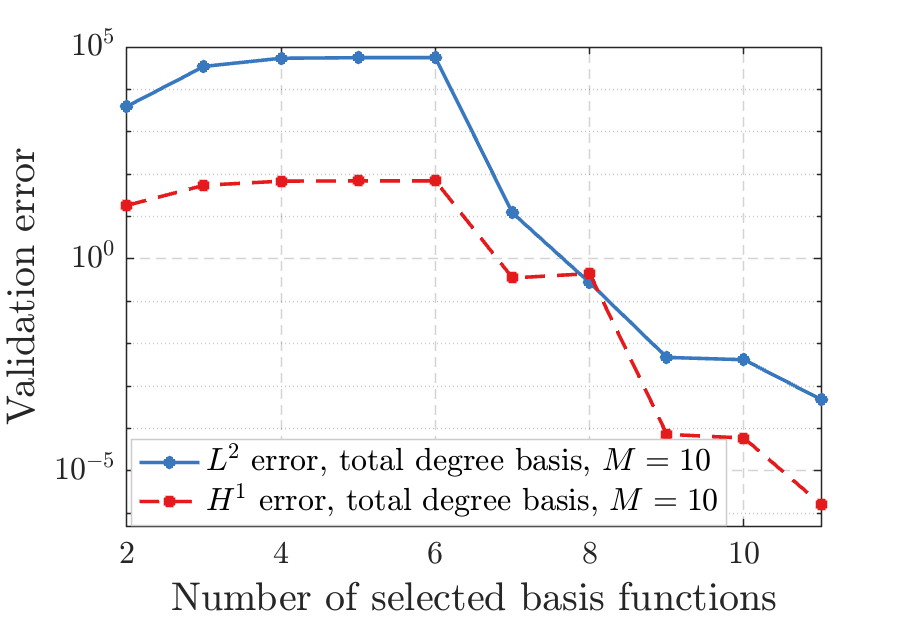} 
  \subcaption{Greedy selection using residual-decay indicator}
\end{subfigure}
\caption{Error decay of Algorithm~\ref{algorithm_1} using the time-dependent Galerkin method for the one-dimensional optimal control problem, comparing random basis selection and greedy selection driven by the residual-decay indicator.}\label{test1d}
\end{figure}

We next consider the one-dimensional optimal control problem from \cite{Kalise_2018}
\begin{equation*}
\min _{u \in L^2(0, T)} \int_0^T \left( \frac{1}{4}\left(y(t)^2 \mathrm{e}^{y(t)}+2 y(t) \mathrm{e}^{y(t)}+4 y(t)^3\right)^2 + u(t)^2 \right) dt
\end{equation*}
subject to 
\begin{equation*}
\dot{y}(t)=u(t),\quad y(0)=x,
\end{equation*}
that is, \eqref{eq:optimal_cost} with $t=0$ and no terminal penalty, $\ell_T\equiv 0$.
The corresponding infinite-horizon stationary HJB equation has the exact solution $V(x)=x^4+x^2\mathrm{e}^x$ on $\Omega=[-1,1]$; for the finite-horizon computation below, $T$ is chosen large enough that $V(0,x)$ is close to this stationary value function on the computational domain. According to Definition~\ref{def:HJB}, given the current approximation $\hat{V}^{(k)}$, the selection indicator for each candidate basis function $\phi_j(x)$ with index $j \in \Lambda^{(k)}$ is
\begin{equation}
e_{j,\mathrm{HJB}}^{(k)}=\frac{\left|\int_{-1}^1 R\left(\hat{V}^{(k)}\right)(x) a_k(x) \phi'_j(x) d x\right|}{\left(\int_{-1}^1\left|a_k(x) \phi'_j(x)\right|^2 d x\right)^{1 / 2}},\quad a_k(x)=-\frac{1}{2} (\hat{V}^{(k)})'(x),
\label{eq:1d_hjb_indicator}
\end{equation}
with residual
\begin{equation*}
R\left(\hat{V}^{(k)}\right)(x)=\frac{1}{4}\left(\left(x^2 \mathrm{e}^x+2 x \mathrm{e}^x+4 x^3\right)^2 - \left((\hat{V}^{(k)})'(x)\right)^2\right).
\end{equation*}
By taking $T=10$ and considering a one-dimensional Legendre polynomial basis truncated at degree 10, we obtain the results reported in Figure \ref{test1d}, where we observe that the accuracy  of  the greedy strategy improves rapidly once the basis size reaches six, while, by contrast, random selection attains the same accuracy only in the final stage.

\subsection{Van der Pol oscillator}
\label{sec:test2d}
We now consider the optimal control of the Van der Pol oscillator:
\begin{equation*}
\min _{u \in L^2(0, T)} \int_0^T \left( y_1(t)^2+y_2(t)^2+\beta\, u(t)^2 \right) dt
\end{equation*}
subject to
\begin{equation*}
\begin{cases}
\dot{y}_1(t)=y_2(t), \\
\dot{y}_2(t)=-y_1(t)+y_2(t)\left(1-y_1(t)^2\right)+u(t), \\
\left(y_1(0), y_2(0)\right)=\left(x_1, x_2\right),
\end{cases}
\end{equation*}
where we set $\beta=0.1$, $T=3$, $\ell_T\equiv 0$, and the spatial domain $\Omega=[-3,3]^2$.

\paragraph{Model-based setting.}
According to Definition~\ref{def:HJB}, given the current approximation $\hat{V}^{(k)}$, for each candidate basis function $\phi_j(\mathbf{x})$ with index $j \in \Lambda^{(k)}$, the selection indicator $e_{j,\mathrm{HJB}}^{(k)}$ is defined with
\begin{equation*}
\begin{aligned}
\mathbf{a}_k(\mathbf{x})&=\binom{x_2}{-x_1+x_2\left(1-x_1^2\right)-\frac{1}{2 \beta} \partial_{x_2} \hat{V}^{(k)}(\mathbf{x})},\\
R\left(\hat{V}^{(k)}\right)(\mathbf{x})&=-\frac{1}{4\beta} \left(\partial_{x_2} \hat{V}^{(k)}(\mathbf{x})\right)^2+\partial_{x_1} \hat{V}^{(k)}(\mathbf{x}) x_2\\
&+\partial_{x_2} \hat{V}^{(k)}(\mathbf{x})\left(-x_1+x_2\left(1-x_1^2\right)\right)+x_1^2+x_2^2.
\end{aligned}
\end{equation*}
We compare the performance of the time-dependent Galerkin method with two basis selection strategies: random selection and greedy selection using the residual-decay indicator. In both cases, we consider either a hyperbolic cross basis with index $\mathfrak{J}=4$ (with a total of $P=52$ basis functions) or a total degree basis with maximum degree $M=8$ (with $P=45$ basis functions), both generated by Legendre polynomials. The initial basis function is selected as the constant basis. The results for random basis selection are averaged over 100 Monte Carlo simulations, while the greedy selection is deterministic when the initial basis is fixed. The validation error is computed over 1000 validation data points.  

To study the error decay with respect to the number of selected basis functions, we set $k_{\max}=P$, which is the size of the whole candidate basis set.
\begin{figure}[ht]
\centering
\captionsetup[subfigure]{aboveskip=-1pt,belowskip=-1pt}
  \begin{subfigure}[t]{0.48\textwidth}
\includegraphics[width=\textwidth]{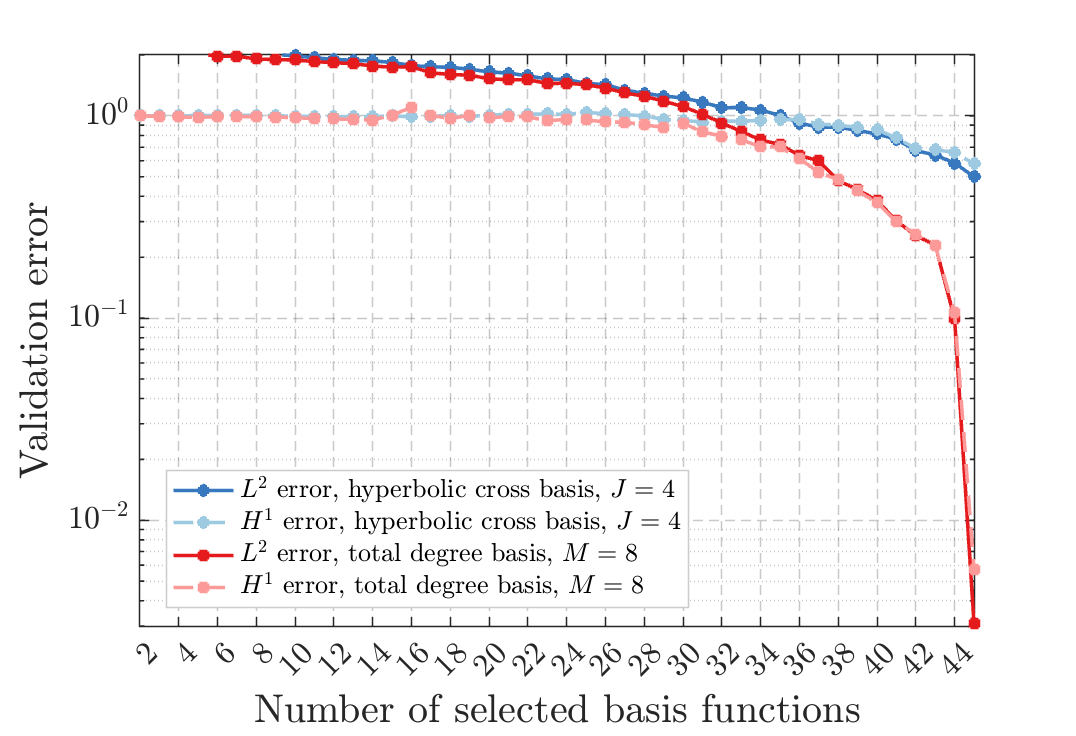}
  \subcaption{Random basis selection}
  \end{subfigure}
\begin{subfigure}[t]
  {0.48\textwidth}
\includegraphics[width=\textwidth]{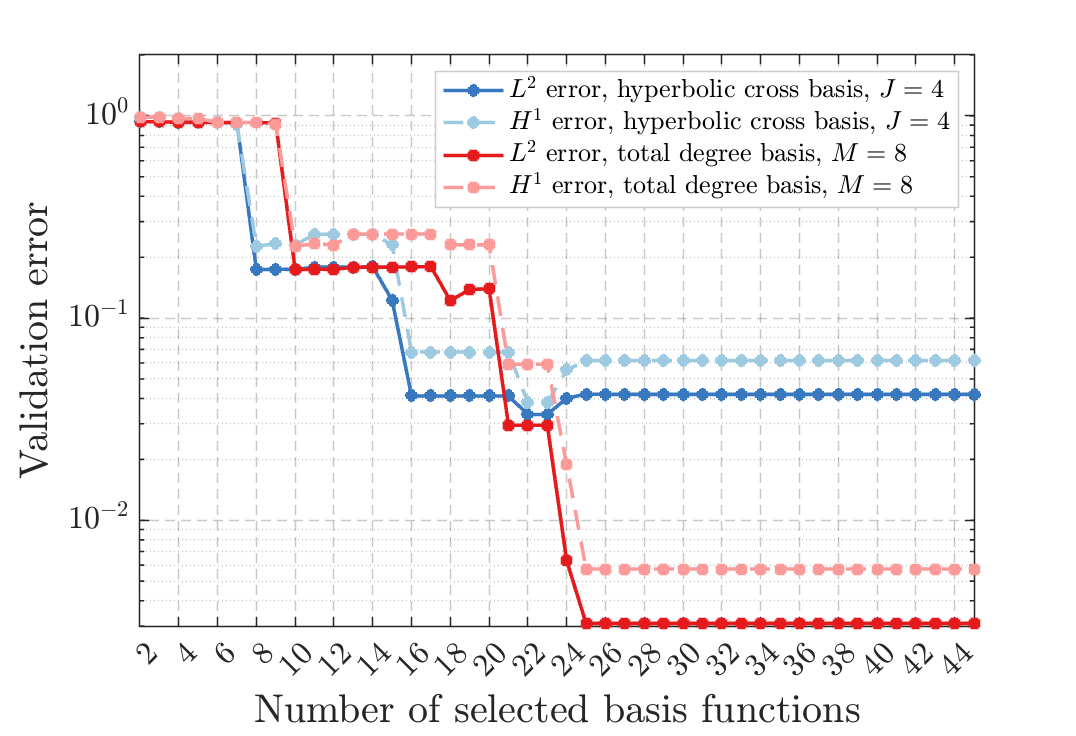} 
  \subcaption{Greedy selection using residual-decay indicator}
\end{subfigure}
\caption{Error decay of Algorithm~\ref{algorithm_1} (model-based basis selection, no data) using the time-dependent Galerkin method for the 2-dimensional Van der Pol oscillator optimal control problem, comparing random and adaptive basis selection strategies.}
\label{fig:galerkin_2d}
\end{figure}
From Figure~\ref{fig:galerkin_2d}, we observe that greedy selection consistently outperforms random selection. Notably, greedy selection achieves the same accuracy level with substantially fewer basis functions, thus improving efficiency by preferentially selecting the most informative basis elements. With the hyperbolic cross basis, greedy selection reaches the accuracy of the full basis with only $16$ of the $52$ candidates, and with the total degree basis with $25$ of the $45$ candidates.

\paragraph{Data-driven setting.}
We now solve the same optimal control problem using the supervised learning framework. Here, training data are used both to fit the model and to compute the supervised loss guiding basis selection. We compare several strategies for jointly selecting data samples and basis functions.

\begin{figure}[h]
    \centering
    \captionsetup[subfigure]{aboveskip=-1pt,belowskip=-1pt}
\begin{subfigure}[t]
  {0.45\textwidth}  \includegraphics[width=\linewidth]{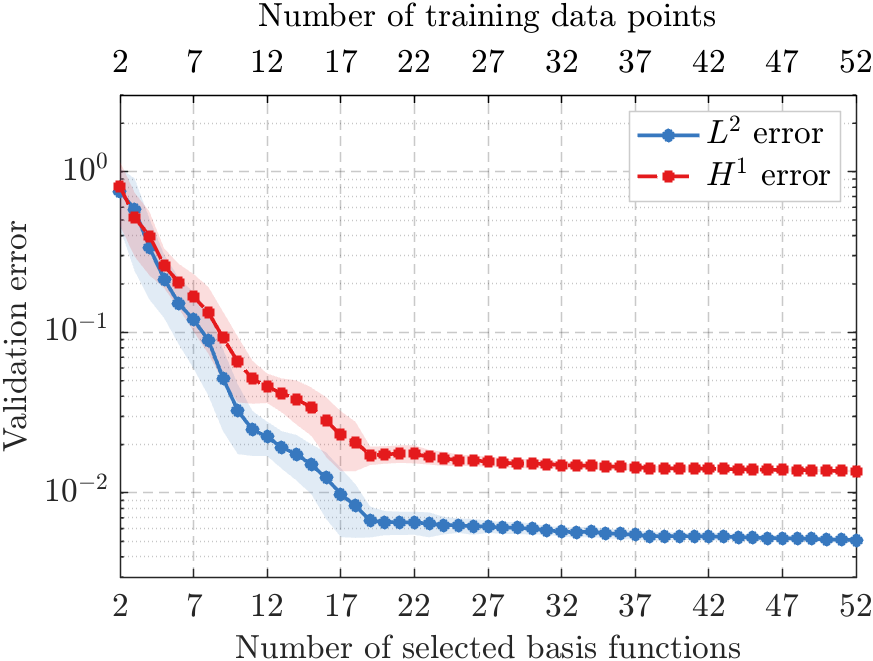}
    \subcaption{Basis selection + data selection}
\end{subfigure}
\begin{subfigure}[t]
  {0.45\textwidth}  \includegraphics[width=\linewidth]{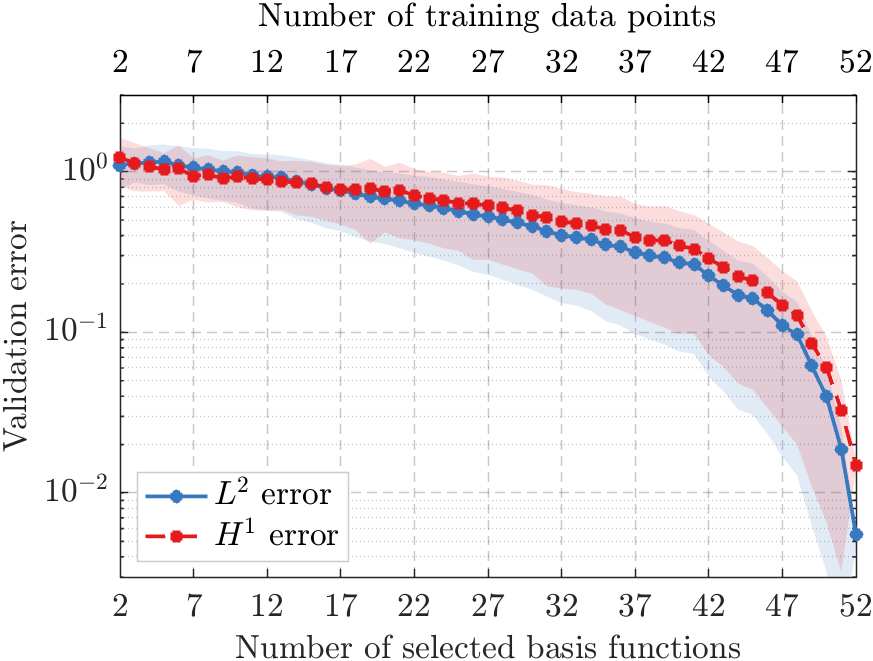}
    \subcaption{Random basis + data selection}
\end{subfigure}
\begin{subfigure}[t]
  {0.45\textwidth}
\includegraphics[width=\linewidth]{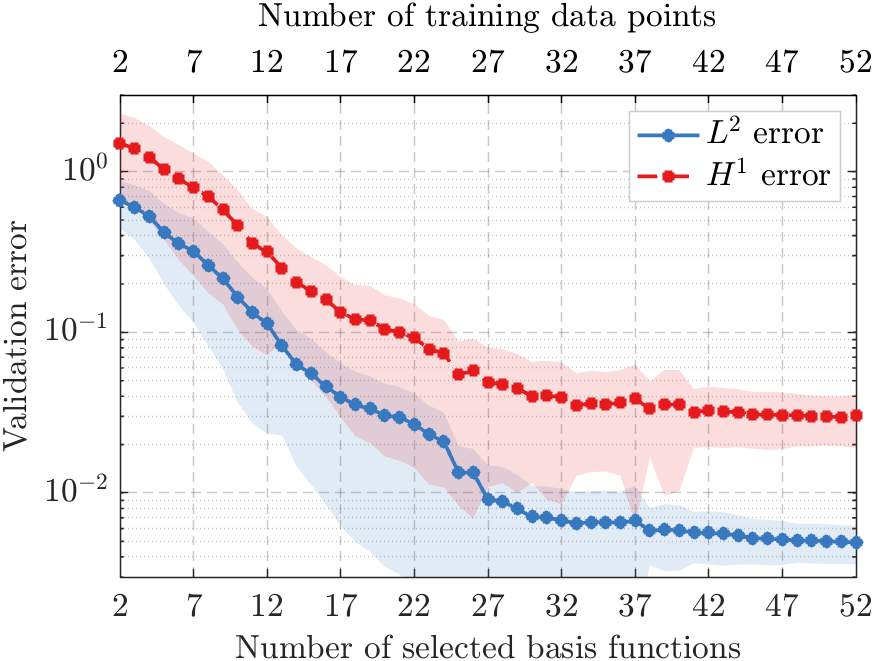}
    \subcaption{Basis selection + random data}
\end{subfigure}
\begin{subfigure}[t]
  {0.45\textwidth}  \includegraphics[width=\linewidth]{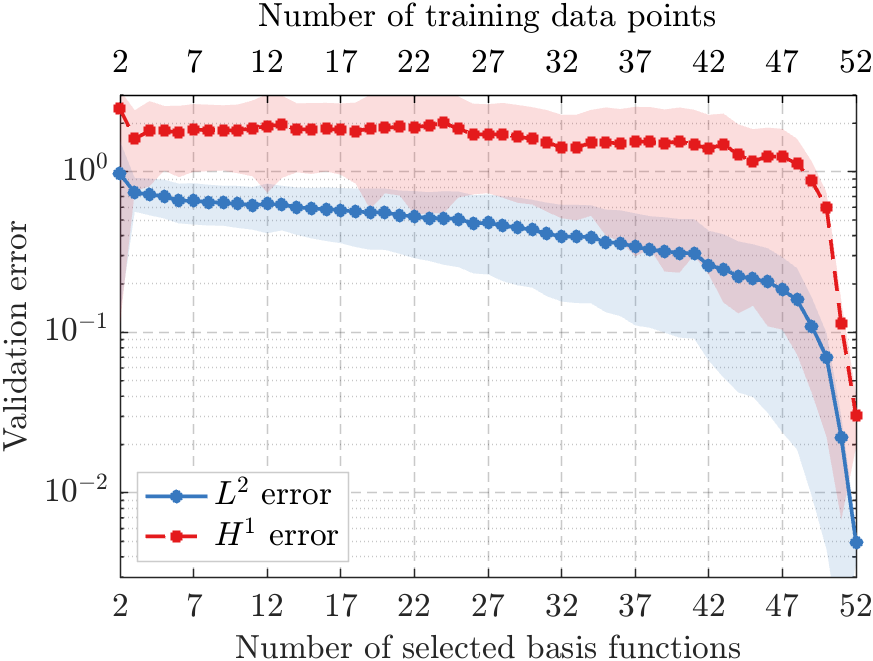}
    \subcaption{Random basis + random data}
\end{subfigure}
\caption{Comparison of adaptive basis and data selection strategies using a hyperbolic cross basis with index $\mathfrak{J}=4$, generated by Legendre polynomials, $k_{\max}=P=52$. Basis functions are selected using the supervised-loss basis indicator. LASSO regression with $\lambda=10^{-4}$ is applied. Results are averaged over 100 independent runs. The shaded regions indicate one standard deviation around the mean validation error curves on a logarithmic scale.}
\label{fig:sp_2d}
\end{figure}

All experiments are initialized with a single randomly selected basis function and data point. In Figure~\ref{fig:sp_2d}, the bottom $x$-axis indicates the number of selected basis functions, while the top $x$-axis shows the number of selected training data points. %For simplicity, we disable early stopping and keep these counts equal at each iteration.
We also take $\mu$ to be the uniform empirical measure over the sample set.

The proposed adaptive basis and data selection strategy is particularly advantageous in the following scenarios:
\begin{itemize}
    \item[i)] \textbf{Expensive data generation.} When obtaining new data points is costly (e.g., via high-fidelity simulations or experiments), it is natural to compare greedy and random basis selection under random data selection, as in Figure~\ref{fig:sp_2d}(c)--(d). This setup reflects a realistic scenario in which data are generated incrementally until sufficiently informative. In this case, adaptive basis selection reaches a prescribed accuracy threshold with fewer data points than random selection, thereby improving data efficiency.
    \item[ii)] \textbf{Representative data selection.} When a large dataset is already available, computational constraints may render it impractical to use the entire set for training or basis updates. An adaptive selection procedure instead identifies a representative subset of data points, together with the most informative basis functions. As shown in Figures~\ref{fig:sp_2d}(a)--(b), adaptive data selection alone provides limited improvement, whereas adaptive basis selection has a significantly greater impact on overall performance, highlighting its critical role in the approximation process.
\end{itemize}

\begin{figure}[h]
    \centering
\includegraphics[width=0.5\linewidth]{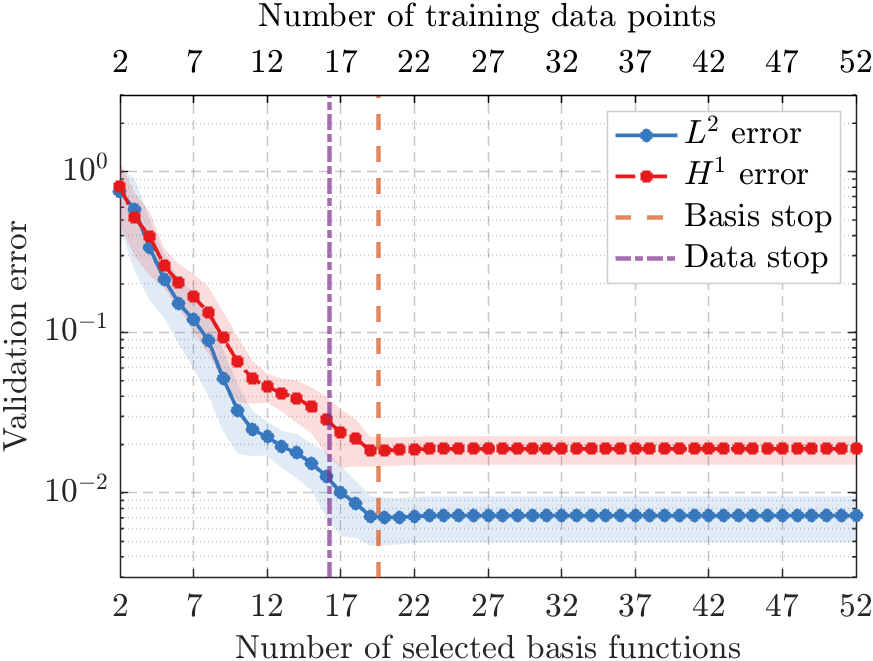}
    \caption{Error decay of Algorithm~\ref{algorithm_2} for the 2-dimensional Van der Pol oscillator optimal control problem with early-stop criteria $\varepsilon_{\mathrm{basis}}=0.05$ and $\varepsilon_{\mathrm{data}}=0.1$. Once the maximum absolute value of indicator/score falls below $\varepsilon$, no further basis functions or data points are added. The vertical lines marks the average stopping point, occurring at 19.54 basis functions and 16.23 data points. Basis functions are selected using the supervised-loss basis indicator.}
    \label{fig:sp_2d_stop}
\end{figure}

\noindent Based on these results, we next activate the early stopping criterion in the adaptive selection procedure, further reducing computational cost while maintaining the desired accuracy. The corresponding results are shown in Figure~\ref{fig:sp_2d_stop}. With $\varepsilon_{\mathrm{basis}}=0.05$ and $\varepsilon_{\mathrm{data}}=0.1$,
Algorithm~\ref{algorithm_2} stops on average after $19$ basis functions and $16$
training samples, that is, only a small subset of the candidate dictionary, while
attaining a validation error comparable to that obtained with the full candidate set.

\subsection{Rigid body model of a satellite}
\label{sec:test6d}
 We now study a six-state rigid body model of a satellite \cite{kang2015causality,kang2017mitigating}.  The state of the satellite is written as $$\mathbf y(t)=[y_1\;\; y_2\;\; y_3\;\; y_4\;\; y_5\;\; y_6]^\top,\qquad
\mathbf y_r:=[y_1\;\; y_2\;\; y_3]^\top,\quad
\mathbf y_\omega:=[y_4\;\; y_5\;\; y_6]^\top,$$
 where $\mathbf y_r$ denotes attitude (rotation angles) and $\mathbf y_\omega$ the  angular velocity in the body-frame. 
The control input is a torque $\mathbf{u}= \left(u_1\;\; u_2\;\; u_3\right)^{\top}$. The optimal control problem is given by
\begin{equation}\label{Test2Eq01}
    \min_{\mathbf{u}\in L^2(0,T;\mathbb{R}^3)} \int^{T}_{0} \left( \alpha_1 \|\mathbf{y}_r(t)\|_2^2 
    + \alpha_2\|\mathbf{y}_\omega(t)\|_2^2 
 + \beta \|\mathbf{u}(t)\|_2^2\right)dt + \alpha_3\left\|\mathbf{y}(T)\right\|_2^2
\end{equation}
subject to
\begin{equation}
\label{Test2Eq02}
\begin{cases}
\dot{\mathbf{y}}_r(t) = \mathbf{E}(\mathbf{y}_r(t))\mathbf{y}_\omega(t), \\
\dot{\mathbf{y}}_\omega(t) = \mathbf{W}\mathbf{S}(\mathbf{y}_\omega(t)) \mathbf{R}(\mathbf{y}_r(t)) \mathbf{h}+\mathbf{W}\mathbf{B}\mathbf{u}(t), \\
\mathbf{y}(0)=\mathbf{x},
\end{cases}
\end{equation}
where $\alpha_1, \alpha_2, \alpha_3,\beta >0$ and the matrices $\mathbf{E}(\mathbf{y}_r)$, $\mathbf{S}(\mathbf{y}_\omega)$, $\mathbf{R}(\mathbf{y}_r) \in \mathbb{R}^{3\times 3}$ are given by
\begin{equation*}
\mathbf{E}(\mathbf{y}_r):=\begin{pmatrix}
1 & \sin y_1 \tan y_2 & \cos y_1 \tan y_2 \\
0 & \cos y_1 & -\sin y_1 \\
0 & \sin y_1 / \cos y_2 & \cos y_1 / \cos y_2
\end{pmatrix},\quad
\mathbf{S}(\mathbf{y}_\omega):= \begin{pmatrix}
0 & y_6 & -y_5 \\
-y_6 & 0 & y_4 \\
y_5 & -y_4 & 0
\end{pmatrix},
\end{equation*}
and
{\scriptsize
\begin{equation*}
\mathbf{R}(\mathbf{y}_r):= \begin{pmatrix}
\cos y_2 \cos y_3 & \cos y_2 \sin y_3 & -\sin y_2 \\
\sin y_1 \sin y_2 \cos y_3-\cos y_1 \sin y_3 & \sin y_1 \sin y_2 \sin y_3+\cos y_1 \cos y_3 & \cos y_2 \sin y_1 \\
\cos y_1 \sin y_2 \cos y_3+\sin y_1 \sin y_3 & \cos y_1 \sin y_2 \sin y_3-\sin y_1 \cos y_3 & \cos y_2 \cos y_1
\end{pmatrix}.
\end{equation*}
}
The constant matrices $\mathbf{W} \in \mathbb{R}^{3\times 3}$, $\mathbf{h} \in \mathbb{R}^3$, and $\mathbf{B} \in \mathbb{R}^{3\times 3}$ are
\begin{equation*}
\mathbf{W} := \begin{pmatrix}
\frac{1}{2} & 0 & 0 \\
0 & \frac{1}{3} & 0 \\
0 & 0 & \frac{1}{4}
\end{pmatrix},\quad \mathbf{h}:=\begin{pmatrix}
1 \\
1 \\
1
\end{pmatrix},\quad \mathbf{B}:= \begin{pmatrix}
1 & \frac{1}{20} & \frac{1}{10} \\
\frac{1}{15} & 1 & \frac{1}{10} \\
\frac{1}{10} & \frac{1}{15} & 1
\end{pmatrix},
\end{equation*}
and $\mathbf{x}:=(x_1, x_2, x_3, x_4, x_5, x_6)^{\top}$ is the initial condition with domain $\Omega = \left[-\tfrac{\pi}{3},\tfrac{\pi}{3}\right]^3 \times \left[-\tfrac{\pi}{4},\tfrac{\pi}{4}\right]^3$. We set $\alpha_1 = \frac{1}{2}$, $\alpha_2 = 5$, $\alpha_3 = \frac{1}{2}$, $\beta = \frac{1}{4}$, and $T=20$. This is an instance of \eqref{eq:optimal_cost} with running cost $\ell(\mathbf{y})=\alpha_1\|\mathbf{y}_r\|_2^2+\alpha_2\|\mathbf{y}_\omega\|_2^2$ and terminal penalty $\ell_T(\mathbf{y})=\alpha_3\|\mathbf{y}\|_2^2$.  

\begin{figure}[ht]
    \centering
    \captionsetup[subfigure]{aboveskip=4pt,belowskip=4pt}
\begin{subfigure}[t]
  {0.45\textwidth}  \includegraphics[width=\linewidth]{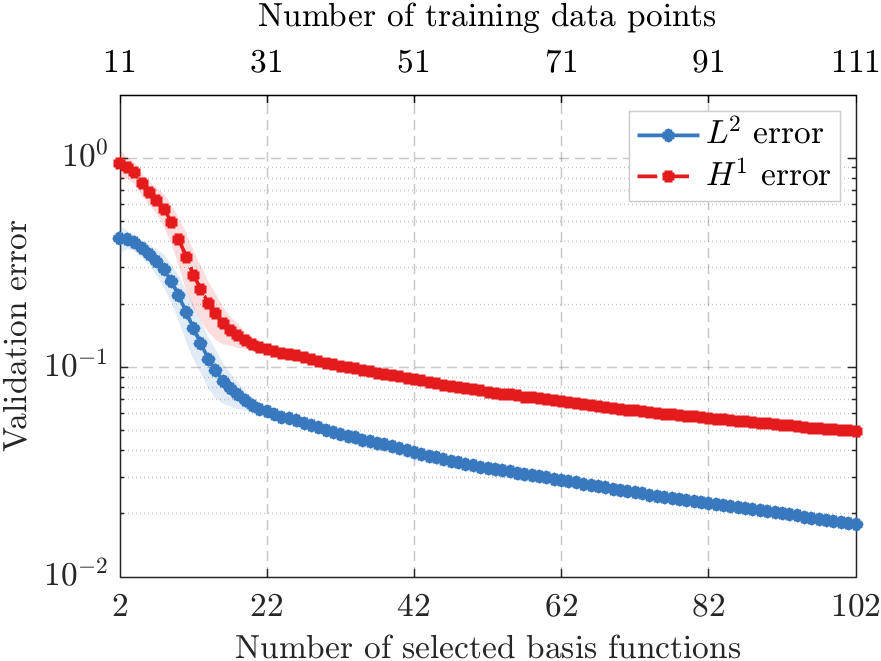}
    \subcaption{Basis selection + data selection}
\end{subfigure}
\begin{subfigure}[t]
  {0.45\textwidth}  \includegraphics[width=\linewidth]{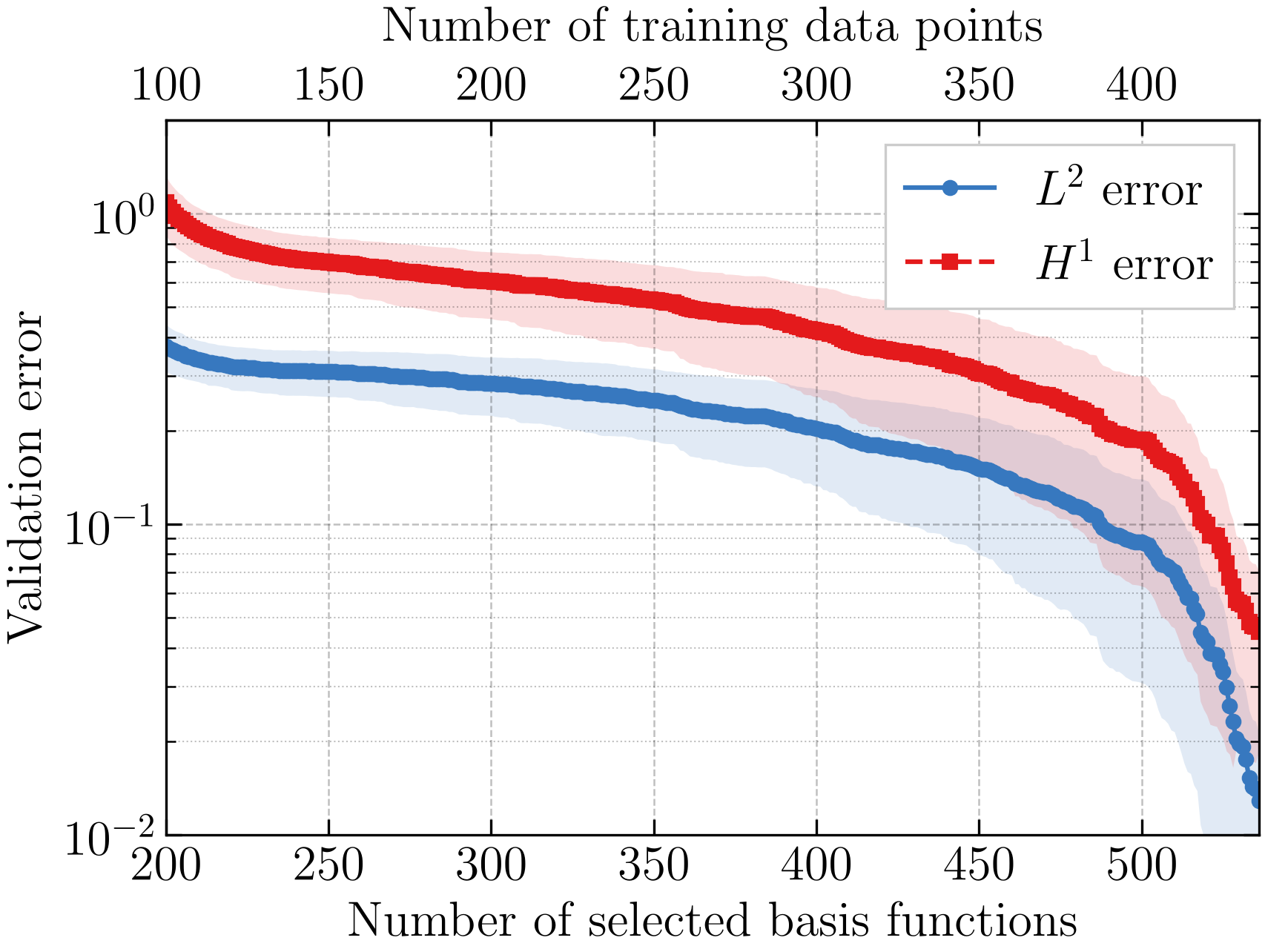}
    \subcaption{Random basis + data selection}
\end{subfigure}
\begin{subfigure}[t]
  {0.51\textwidth}
\includegraphics[width=\linewidth]{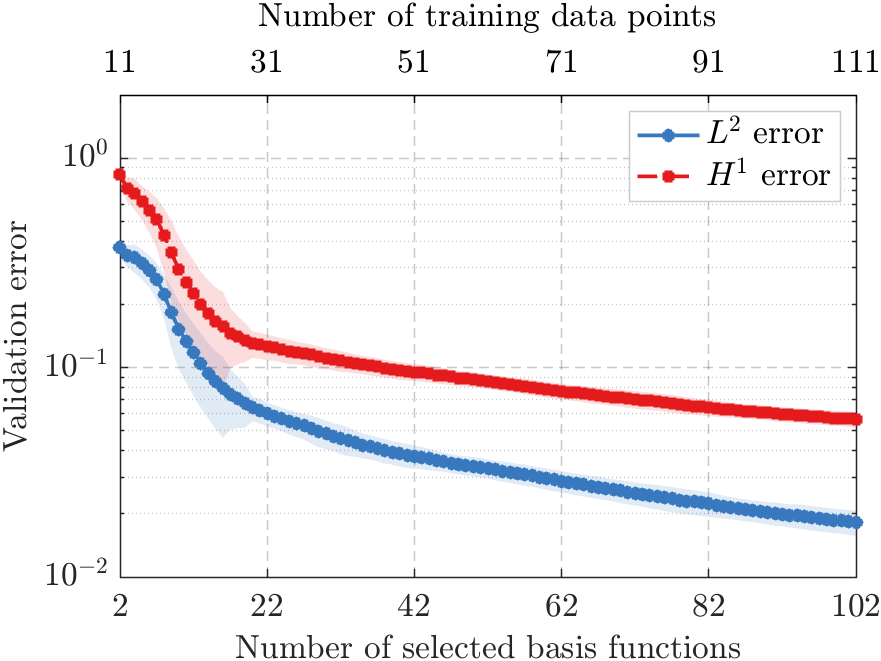}
    \subcaption{Basis selection + random data}
\end{subfigure}
\begin{subfigure}[t]
  {0.45\textwidth}  \includegraphics[width=\linewidth]{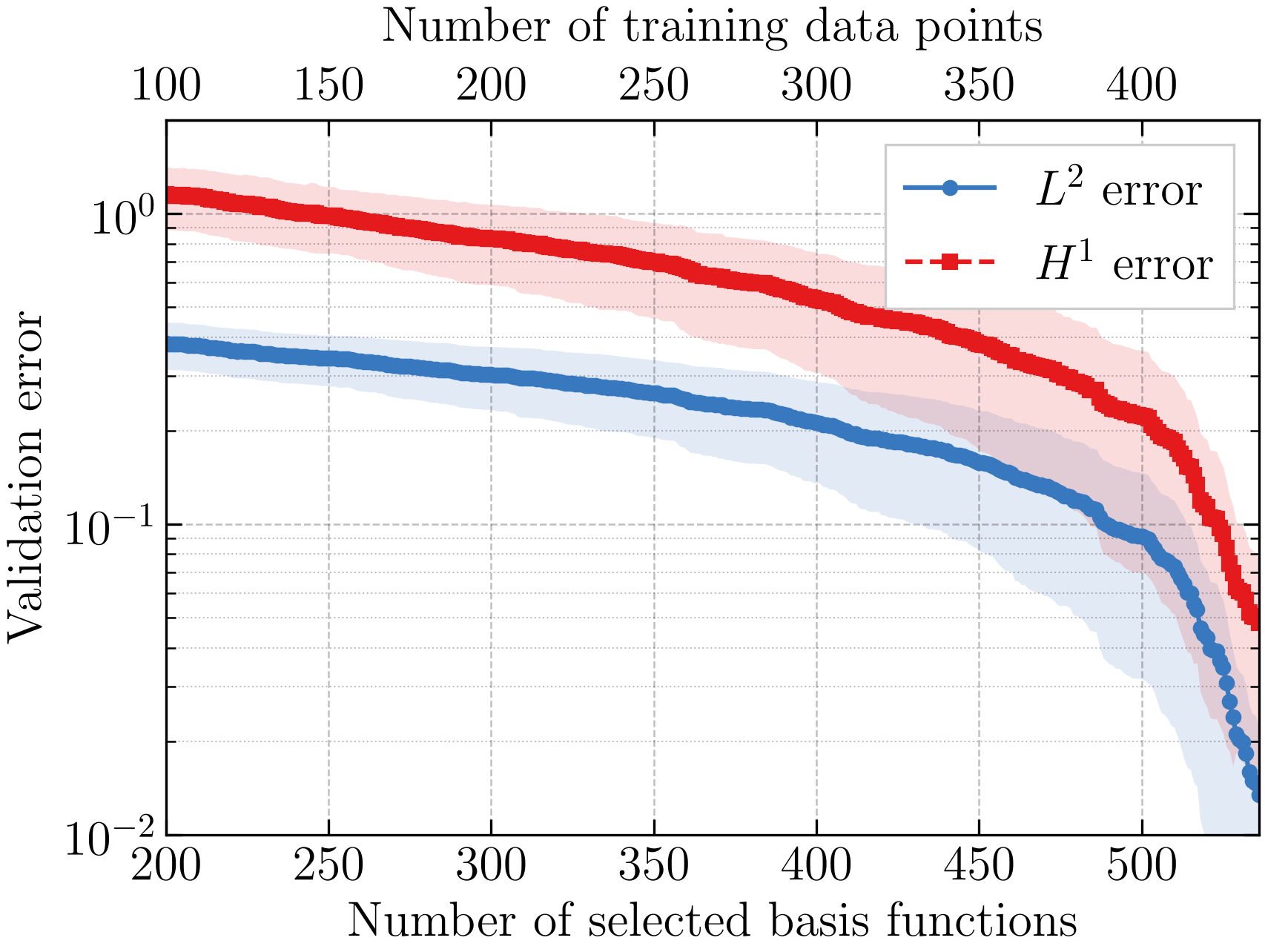}
    \subcaption{Random basis + random data}
\end{subfigure}
\caption{Comparison of adaptive basis and data selection strategies using a hyperbolic cross basis with index $\mathfrak{J}=4$ (the cardinality of the generated basis is 537), generated by Legendre polynomials. Basis functions are selected using the supervised-loss basis indicator. Coefficients are computed by LASSO regression with $\lambda=10^{-4}$. Results are averaged over 100 independent runs, the validation error is computed over 1000 validation points.}
\label{fig:sp_6d}
\end{figure}

For this problem, we focus on the supervised learning framework.  Note that, in the present setting,  the adaptive basis-selection strategies (Figure~\ref{fig:sp_6d}(a) and (c)) are initialized
with a small active set and reach the satisfactory accuracy after 
roughly $10^2$ iterations, while the non-adaptive basis baselines (Figure~\ref{fig:sp_6d}(b) and (d)) attain a
comparable error level only once a substantial fraction of the $P = 537$
candidate functions has been included, and are therefore plotted over the
correspondingly larger range.

\subsection{35-dimensional infinite-horizon optimal control problem}
To further assess the performance of the proposed adaptive basis and data selection
strategy in higher dimensions, we consider an infinite-horizon optimal control problem with a known exact solution. 

We consider the infinite-horizon problem \eqref{eq:infinite_horizon} with $m=d$, $\mathbf{q}(\mathbf{x})=\mathbf{I}_d$, and linear dynamics $\mathbf{f}(\mathbf{x})=(f_1(\mathbf{x}),\dots,f_d(\mathbf{x}))^\top$ given by
$$
f_i(\mathbf{x})=-\kappa_i x_i, \qquad \kappa_i = 0.3+\frac{0.5(i-1)}{d-1},\qquad i=1,\ldots,d.
$$
The computational domain is $\Omega=[-1,1]^d$, and in this experiment we take
$d=35$ and $\beta=1$.  We prescribe the exact value function
$$
V(\mathbf{x}) = 1+ \sum_{i=1}^{6} \rho_i x_i^2 + \sum_{i=1}^{6} \eta_i x_i^3 + 0.08\sum_{i=1}^{5}x_i x_{i+1}, \qquad \mathbf{x}\in\Omega, $$
where $\rho_i = 1.2 - 0.13(i-1)$ and $\eta_i = 0.06\big(1 - 0.11(i-1)\big)$ for $i = 1,\dots,6$, and $\ell(\mathbf{x})$ is chosen accordingly. This test yields a smooth solution and admits a sparse representation in the polynomial basis.

Training and validation datasets are generated by sampling points uniformly from
$\Omega$. We use the hyperbolic cross basis with index $\mathfrak{J}=3$, generated by Legendre polynomials. The adaptive method is initialized with $10$ randomly selected training samples and the constant basis function, and performs $50$ iterations, each adding one training sample and one basis function. After each enrichment step, the coefficients are recomputed on the current active basis and data sets using LASSO regression with regularization parameter $\lambda=10^{-8}$.
For basis selection, we compare the supervised-loss indicator \eqref{eq:supervised_loss_indicator}, with weights $\omega_1=1$ and $\omega_2=1/d$, and the empirical residual-decay indicator \eqref{eq:empiricalPDEResidue}. Each is combined with either the data score \eqref{eq:data_score} or random sampling for data selection. As a baseline, we also fit the full basis on the same randomly selected training samples.

\begin{figure}[!htbp]
\centering
\captionsetup[subfigure]{aboveskip=3pt,belowskip=4pt}
\begin{subfigure}[t]{0.45\textwidth}
  \includegraphics[width=\linewidth]{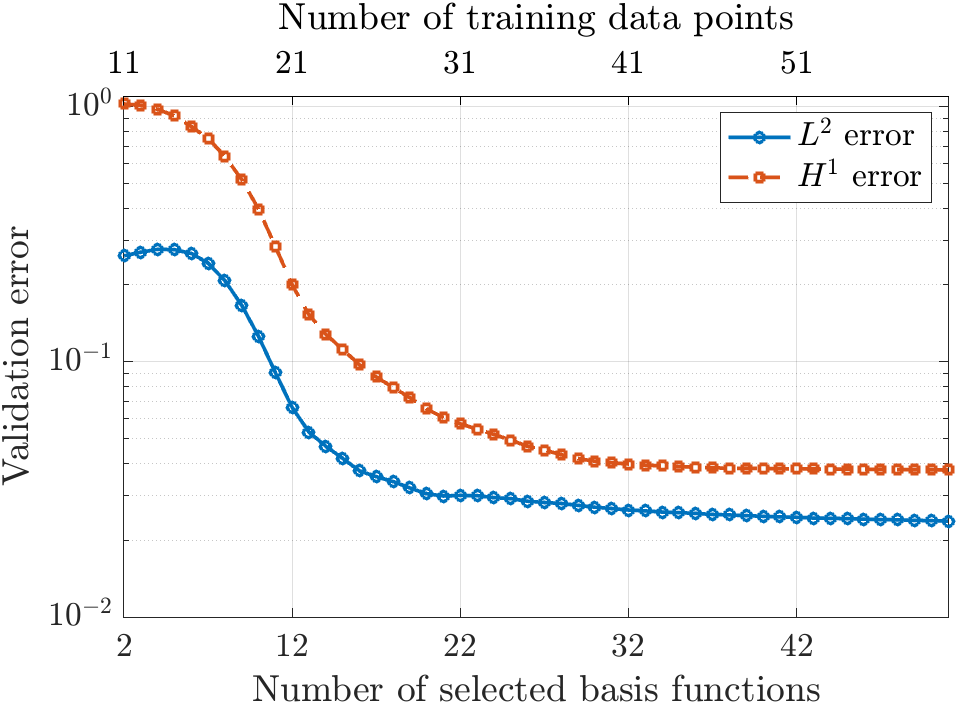}
  \subcaption{$e_{j,\mathrm{SL}}$ + data selection}
\end{subfigure}
\begin{subfigure}[t]{0.45\textwidth}
  \includegraphics[width=\linewidth]{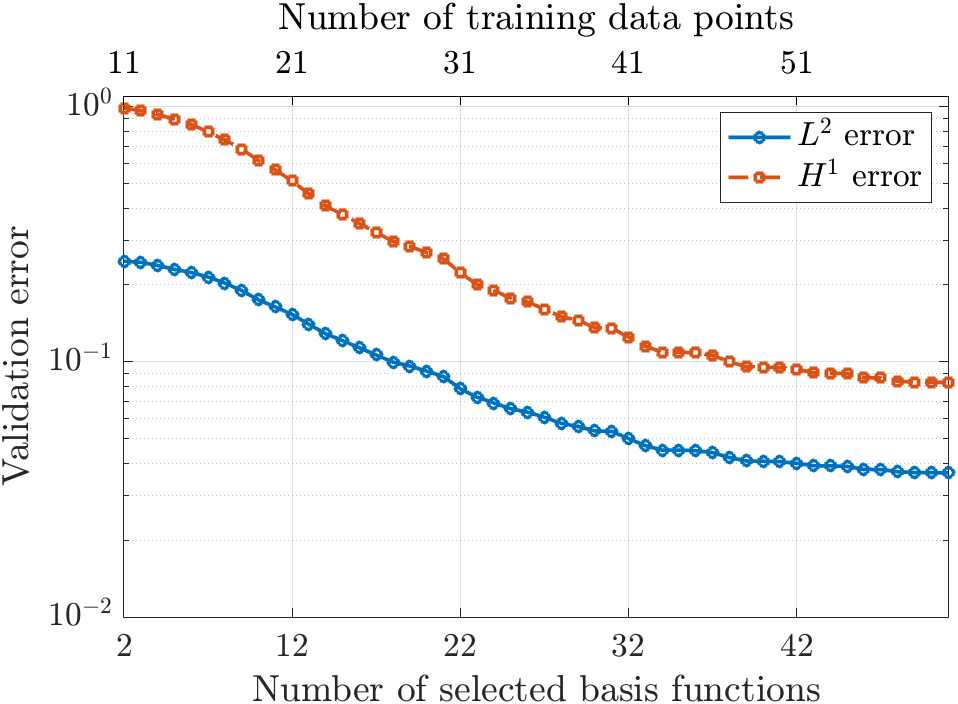}
  \subcaption{$e_{j,\mathrm{SL}}$ + random data}
\end{subfigure}
\begin{subfigure}[t]{0.45\textwidth}
  \includegraphics[width=\linewidth]{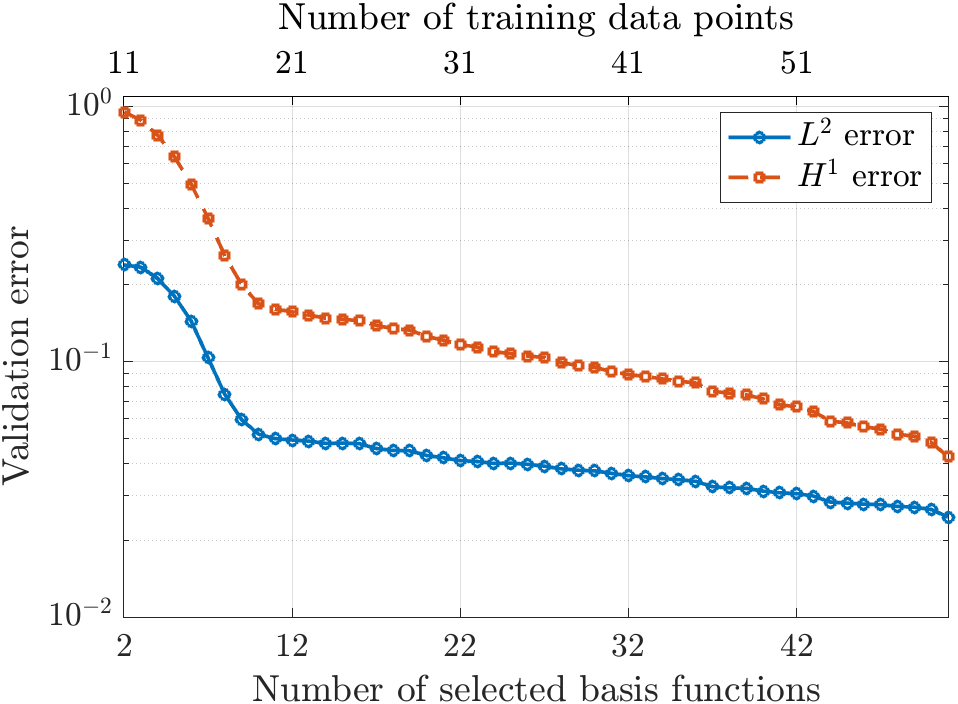}
  \subcaption{$\tilde e_{j,\mathrm{PDE}}$ + data selection}
\end{subfigure}
\begin{subfigure}[t]{0.45\textwidth}
  \includegraphics[width=\linewidth]{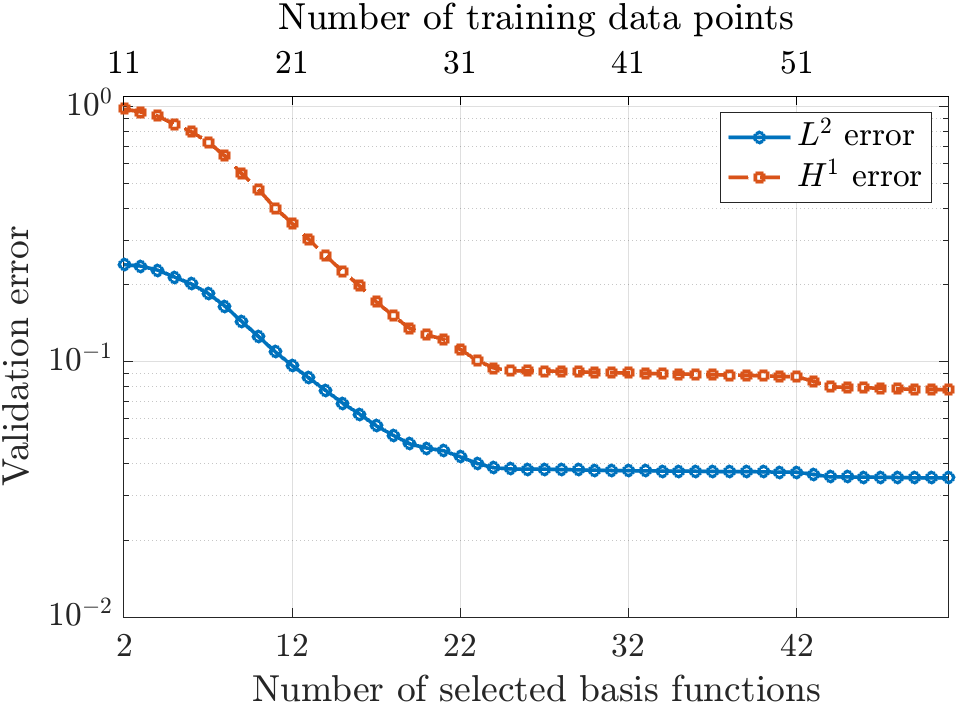}
  \subcaption{$\tilde e_{j,\mathrm{PDE}}$ + random data}
\end{subfigure}
\begin{subfigure}[t]{0.45\textwidth}
  \includegraphics[width=\linewidth]{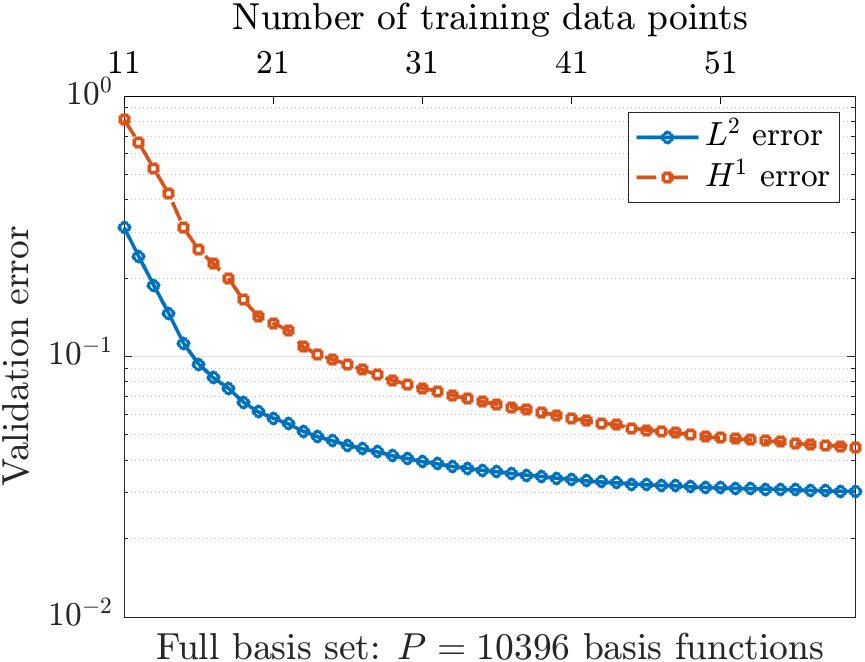}
  \subcaption{Full basis + random data}
\end{subfigure}
\caption{Validation error decay for the stationary HJB problem with $d=35$. In (a)--(d), Algorithm~\ref{algorithm_2} starts from $10$ random training samples and the constant basis function, and adds one sample and one basis function per iteration. Basis functions are selected with the supervised-loss indicator \eqref{eq:supervised_loss_indicator} in (a)--(b), and with the empirical residual-decay indicator \eqref{eq:empiricalPDEResidue} in (c)--(d). Training samples are selected with the data score \eqref{eq:data_score} in (a) and (c), and at random in (b) and (d). Subfigure (e) shows the full-basis baseline, which uses all $10396$ basis functions, with training samples selected randomly from the same initial size.}
\label{fig:35d}
\end{figure}

It can be seen from Figure~\ref{fig:35d} that, when combined with the data score \eqref{eq:data_score}, both basis indicators achieve an error decay comparable to the full-basis regression, which uses the full candidate dictionary of $P=10396$ basis functions. Notably, the adaptive strategy attains similar accuracy with an active subspace of only $51$ basis functions, demonstrating its computational efficiency in high-dimensional settings. In contrast, when the training samples are selected at random, the error decays more slowly for both indicators.

\section{Conclusion}
We have introduced two adaptive strategies for approximating PDE solutions in both model-based and data-driven regimes, with a particular emphasis on high-dimensional nonlinear HJB equations. The Galerkin framework utilizes a closed-form residual-decay indicator to greedily enrich the basis. The supervised learning framework combines greedy basis expansion with adaptive selection of gradient-augmented training samples driven by a supervised loss, focusing computational effort on the most informative data points. Across numerical experiments on benchmark and control problems, both frameworks achieve accurate solutions with substantially fewer basis functions and far fewer data points (or without data).

These proposed residual-decay indicators, derived from a dual-weighted residual argument via single-direction orthogonal projection, effectively guide basis and data selection in high-dimensional settings. Unlike classical goal-oriented DWR estimators, which aim to control the error in a user-specified functional of physical interest, our indicators are not designed as rigorous a posteriori error estimators. Instead, they provide simple and computationally efficient criteria for adaptive basis enrichment. As a consequence, the predicted residual decay does not necessarily coincide with the actual reduction in approximation error, even though it proves highly effective in practice.

Several directions remain for future work. These include developing quantitative stopping criteria for adaptive basis and data selection, and extending the current single-direction enrichment strategy to multi-direction updates, in which batches of basis functions are selected at each iteration, potentially improving both numerical performance and overall efficiency.

\noindent\subsection*{Acknowledgments} Y.H. is supported by a Research Fellowship from the Royal Commission for the Exhibition of 1851. D.K. is supported by the Air Force Office of Scientific Research under award number FA8655-26-1-B011. Y.Z. is supported by the National Science Foundation under DMS award number 2529107.
 
\bibliographystyle{spmpsci}
\bibliography{references}
\end{document}